\documentclass{article}
\usepackage{latexsym,amsfonts,amsmath,amsthm,makeidx}

\let\oldbibliography\thebibliography
\renewcommand{\thebibliography}[1]{%
\oldbibliography{#1}%
\setlength{\itemsep}{0pt}%
}

\makeindex
\newtheorem{definition}{Definition}[section]
\newtheorem{theorem}{Theorem}[section]
\newtheorem{lemma}{Lemma}[section]
\newtheorem{corollary}{Corollary}[section]
\newtheorem{proposition}{Proposition}[section]
\newtheorem{remark}{Remark}[section]
\newcommand{\RN}{\mathbb R^N}
\newcommand{\om}{\Omega}

\newcommand{\iy}{\infty}

\newcommand{\s}{\section}
\newcommand{\dd}{\delta}

\newcommand{\la}{\lambda}

\newcommand{\p}{\phi}

\newcommand{\R}{\mathbb R}
\newcommand{\al}{\alpha}

\newcommand{\bb}{\beta}

\newcommand{\e}{\varepsilon}
\newcommand{\vp}{\varphi}

\newcommand{\bt}{\begin{theorem}}
\newcommand{\et}{\end{theorem}}
\newcommand{\bl}{\begin{lemma}}
\newcommand{\el}{\end{lemma}}
\newcommand{\bd}{\begin{definition}}
\newcommand{\ed}{\end{definition}}
\newcommand{\bc}{\begin{corollary}}
\newcommand{\ec}{\end{corollary}}
\newcommand{\bp}{\begin{proof}}
\newcommand{\ep}{\end{proof}}
\newcommand{\bx}{\begin{example}}
\newcommand{\ex}{\end{example}}
\newcommand{\bi}{\begin{exercise}}
\newcommand{\ei}{\end{exercise}}
\newcommand{\bo}{\begin{prop}}
\newcommand{\eo}{\end{prop}}
\newcommand{\br}{\begin{remark}}
\newcommand{\er}{\end{remark}}
\newcommand{\be}{\begin{equation}}
\newcommand{\ee}{\end{equation}}
\newcommand{\ba}{\begin{align}}
\newcommand{\ea}{\end{align}}
\newcommand{\bn}{\begin{enumerate}}
\newcommand{\en}{\end{enumerate}}
\newcommand{\bg}{\begin{align*}}
\newcommand{\bcs}{\begin{cases}}
\newcommand{\ecs}{\end{cases}}

\newcommand{\intR}[1]{\int_{\RN}#1\,  dx}

\newcommand{\va}{\varphi}
\newcommand{\sg}{\sigma}
\newcommand{\bean}{\begin{eqnarray*}}
\newcommand{\eean}{\end{eqnarray*}}

\numberwithin{equation}{section}

\begin{document}

\title{\bf Positive Least Energy Solutions and Phase Separation for Coupled
Schr\"{o}dinger Equations with Critical Exponent: Higher Dimensional Case\thanks{Supported by NSFC (11025106).
E-mail:
chenzhijie1987@sina.com(Chen);\quad
\quad wzou@math.tsinghua.edu.cn
(Zou)}}
\date{}
\author{{\bf Zhijie Chen,
  Wenming Zou}\\
\footnotesize {\it  Department of Mathematical Sciences, Tsinghua
University,}\\
\footnotesize {\it Beijing 100084, China}}

\maketitle
\begin{center}
\begin{minipage}{120mm}
\begin{center}{\bf Abstract}\end{center}

We study the following nonlinear Schr\"{o}dinger system which is related to Bose-Einstein
condensate:
\begin{displaymath}
\begin{cases}-\Delta u +\la_1 u =
\mu_1 u^{2^\ast-1}+\beta u^{\frac{2^\ast}{2}-1}v^{\frac{2^\ast}{2}}, \quad x\in \Omega,\\
-\Delta v +\la_2 v =\mu_2 v^{2^\ast-1}+\beta v^{\frac{2^\ast}{2}-1} u^{\frac{2^\ast}{2}},  \quad   x\in
\om,\\
u\ge 0, v\ge 0 \,\,\hbox{in $\om$},\quad
u=v=0  \,\,\hbox{on $\partial\om$}.\end{cases}\end{displaymath}
Here $\om\subset \R^N$ is a smooth bounded domain,
$2^\ast:=\frac{2N}{N-2}$ is the Sobolev critical exponent, $-\la_1(\om)<\la_1,\la_2<0$, $\mu_1,\mu_2>0$ and $\beta\neq 0$,
where $\lambda_1(\om)$
is the first eigenvalue of $-\Delta$ with the Dirichlet
boundary condition. When $\bb=0$, this is just the well-known Brezis-Nirenberg problem.
The special case $N=4$ was studied by the authors in (Arch. Ration. Mech. Anal. 205: 515-551, 2012). In this paper we consider {\it the higher dimensional case $N\ge 5$}. It is interesting that
we can prove the existence of a positive least energy solution $(u_\bb, v_\bb)$ {\it for any $\beta\neq 0$ } (which can not hold in the special case $N=4$).
We also study the limit behavior of $(u_\bb, v_\bb)$ as $\beta\to -\infty$ and
phase separation is expected. In particular,
$u_\bb-v_\bb$ will converge to {\it sign-changing solutions} of the Brezis-Nirenberg problem,
provided $N\ge 6$. In case $\la_1=\la_2$,
the classification of the least energy solutions is also studied.
It turns out that some quite different phenomena appear comparing to the special case $N=4$.\\
\end{minipage}
\end{center}

 \s{Introduction}

In this paper we consider the following coupled nonlinear
Schr\"{o}dinger equations:
\be\label{eq2}
\begin{cases}-\Delta u +\la_1 u =
\mu_1 u^{2p-1}+\beta u^{p-1}v^{p}, \quad x\in \Omega,\\
-\Delta v +\la_2 v =\mu_2 v^{2p-1}+\beta v^{p-1} u^{p},  \quad   x\in
\om,\\
u\ge 0, v\ge 0 \,\,\hbox{in $\om$},\quad
u=v=0  \,\,\hbox{on $\partial\om$},\end{cases}\ee
where $\om=\RN$ or $\om\subset\RN$ is a smooth bounded domain, $p>1$ and $p\le 2^\ast/2$ if $N\ge 3$,
$\mu_1,\mu_2 >0$ and $\beta\neq 0$ is a coupling constant.
In the case $p=2$, the cubic system (\ref{eq2}) appears in many physical problems,
especially in nonlinear optics and Bose-Einstein condensation. We refer for this to \cite{AA, EGBB,F, KL}, which
also contain information about the physical relevance of non-cubic nonlinearities.
In the subcritical case $p<2^\ast/2$, the existence of solutions have
received great interest recently, see \cite{AC2, BDW, DWW, LW1, LW,MMP, S, WW2} and references therein.

All the papers mentioned above deal with the subcritical case. To the
best of our knowledge, there are no existence results for (\ref{eq2}) in the critical case $2p=2^\ast$ in
the literature. Critical exponent problems, which have received great attention in the past thirty years
since the cerebrated work by Brezis and Nirenberg \cite{BN},
are very interesting in view of mathematics. Recently, the authors studied
the special critical case $p=2$ and $N=4$ in \cite{CZ}.

In this paper, we study the existence and properties of least energy solutions to (\ref{eq2})
in {\it the higher dimensional case}. In the sequel we assume that
\be\label{fc} N\ge 5 \quad\hbox{and}\quad 2p=2^\ast.\ee
It turns out that different phenomena happen comparing to the special case $N=4$ (\cite{CZ}),
see Remarks \ref{rmk}, \ref{remark}, \ref{remark3} and \ref{rmk1} below.
If $\om=\R^N$ and $(u, v)$ is any
a solution of (\ref{eq2}), then by the Pohozaev identity, it is easy to get that
$\intR{\la_1 u^2+\la_2 v^2}=0,$ so $(u, v)\equiv(0, 0)$ if $\la_1 \la_2>0$. Therefore,
in the sequel we assume that {\it $\om\subset \R^N$ is a smooth bounded domain}
. We call a solution $(u, v)$ nontrivial
if both $u\not\equiv 0$ and $v\not\equiv 0$, a solution $(u, v)$ semi-trivial if $(u, v)$ is type
of $(u, 0)$ or $(0, v)$. We are concerned with nontrivial solutions of (\ref{eq2}).

Let $\lambda_1(\om)$ be the first eigenvalue of $-\Delta$ with the Dirichlet
boundary condition. Clearly (\ref{eq2}) has semi-trivial
solutions $(u_{\mu_1}, 0)$ and $(0, u_{\mu_2})$, where
$u_{\mu_i}$ is a positive least energy solution of the well-known Brezis-Nirenberg problem
\be\label{Brezis-Nirenberg}-\Delta u +\la_i u=\mu_i |u|^{2^\ast-2}u,\,\, u\in H^1_0(\om),\ee
if $-\la_1(\om)<\la_1, \la_2<0$ (see \cite{BN}). Hence, system (\ref{eq2}) can be seen as
a critically coupled Brezis-Nirenberg problem. As we will see in Theorems \ref{th3}-\ref{sign-changing}, system (\ref{eq2}) is closely related to the Brezis-Nirenberg problem (\ref{Brezis-Nirenberg}). The Brezis-Nirenberg problem (\ref{Brezis-Nirenberg}) has
been studied intensively,
and we refer the readers to
\cite{ CFS, CSS,CZ1, DS, SZ} and references therein.

Denote $H:=H^1_0(\om)\times H^1_0(\om)$.
 It is well known that solutions
of (\ref{eq2}) correspond to the critical points of $C^1$ functional $E: H\to \R$
given by
{\allowdisplaybreaks
\begin{align}\label{eq3}
E(u, v)=&\frac{1}{2}\int_{\om}(|\nabla u|^2+\la_1 u^2)+\frac{1}{2}\int_{\om}(|\nabla v|^2+\la_2 v^2)\nonumber\\
&-\frac{1}{2p}\int_{\om}(\mu_1 |u|^{2p}+2\beta |u|^{p}|v|^{p} +\mu_2 |v|^{2p}).
\end{align}
}%

We say a solution $(u, v)$ of (\ref{eq2}) is {\it a least energy solution}, if $(u, v)$ is nontrivial and $E(u, v)\le E(\varphi,\psi)$ for
any other nontrivial solution $(\varphi,\psi)$ of (\ref{eq2}). As in \cite{LW1}, we define a Nehari type manifold
{\allowdisplaybreaks
\begin{align*}\mathcal{M}=\left\{(u, v)\in H\,\,:\,\,u\not\equiv 0,\, v\not\equiv 0,\, E'(u, v)(u, 0)=E'(u, v)(0, v)=0\right\}.\end{align*}
}%
Then any nontrivial solutions of (\ref{eq2}) belong to $\mathcal{M}$.
Similarly as \cite{CZ}, it is trivial to see that $\mathcal{M}\neq \emptyset$. Define the least energy
\be\label{eq10}B:=\inf_{(u, v)\in \mathcal{M}} E(u, v)=\inf_{(u, v)\in \mathcal{M}}
\frac{1}{N}\int_{\om}(|\nabla u|^2+\la_1u^2+|\nabla v|^2+\la_2 v^2)\,dx.\ee

First we consider the symmetric case $-\la_1(\om)<\la_1=\la_2=\la<0$. By \cite{BN} the Brezis-Nirenberg problem
\be\label{BN}-\Delta u +\la u= |u|^{2^\ast-2}u,\,\, u\in H^1_0(\om)\ee
has a positive least energy solution $\omega$ with energy
\be\label{func}B_1 :=\frac{1}{N}\int_{\om}(|\nabla \omega|^2+\la \omega^2)\,dx=\frac{1}{N}\int_{\om}\omega^{2^\ast}\,dx.\ee
Moreover,
\be\label{func1}\int_{\om}(|\nabla u|^2+\la u^2)\,dx\ge (N B_1)^{2/N}\left(\int_{\om}|u|^{2^\ast}\,dx\right)^{2/2^\ast}, \quad\forall\,\, u\in H^1_0(\om).\ee
Consider the following nonlinear problem ($p=\frac{N}{N-2}<2$ since $N\ge 5$)
\be\label{eq8}
\begin{cases}\mu_1 k^{p-1}+\beta  k^{p/2-1}l^{p/2}  =
1,\\
\beta k^{p/2} l^{p/2-1}+\mu_2 l^{p-1}  =1,\\
k>0,\,\,\, l>0.\end{cases}\ee
We will prove in Lemma \ref{lem} that there exists $(k_0, l_0)$ such that
\be\label{eqqqq} \hbox{$(k_0, l_0)$ satisfies (\ref{eq8}) and}\,\,\,   k_0=\min\{k\,: \,\hbox{$(k, l)$ is a solution of (\ref{eq8})}\}.\ee
Our first result deals with the symmetric case $\la_1=\la_2$.

\bt\label{th0} Assume that $-\la_1(\om)<\la_1= \la_2=\la<0$. Let $(k_0, l_0)$ in (\ref{eqqqq}).
Then for any $\bb>0$,  $(\sqrt{k_0}\omega,\sqrt{l_0}\omega)$ is a positive solution of (\ref{eq2}).
Moreover, if $\bb\ge \frac{2}{N-2}\max\{\mu_1, \mu_2\}$, then $E(\sqrt{k_0}\omega, \sqrt{l_0}\omega)=B$, that is,
 $(\sqrt{k_0}\omega, \sqrt{l_0}\omega)$ is a positive least energy solution of (\ref{eq2}).
\et

\br\label{rmk} (1) In the special case $N=4$ and $2p=2^\ast$, \cite[Theorem 1.1]{CZ} said that (\ref{eq2}) has no nontrivial nonnegative solution if
$\beta\in [\min\{\mu_1,\mu_2\}, \max\{\mu_1,\mu_2\}] $ and $\mu_1\neq \mu_2$.
Therefore, the general case $N\ge 5$ is quite different from
the case $N=4$. As we will see in Section 2, the proof of Theorem \ref{th0} is much more delicate than
 \cite{CZ} because of the nonlinearity of problem (\ref{eq8}).

(2) Similarly as in \cite[Remark 1.1]{CZ}, we can prove that, if $\om$ is starshaped, the assumption $-\la_1(\om)<\la<0$
in Theorem \ref{th0} is optimal.\er

Our second result deals with the classification of the least energy solutions.

\bt\label{unique1} Let assumptions in Theorem \ref{th0} hold. There exists $\bb_0\ge \frac{2}{N-2}\max\{\mu_1, \mu_2\}$
determined by $(\mu_1, \mu_2)$, and assume that $\bb>\bb_0$. Let $(u, v)$ be any a positive least energy solution
of (\ref{eq2}), then $(u, v)=(\sqrt{k_0}U, \sqrt{l_0}U)$, where $U$ is a positive least energy solution of (\ref{BN}). In particular, the positive least energy
solution of (\ref{eq2}) is unique if $\om\subset \RN$ is a ball.\et

\br (1) We can give a precise definition of $\bb_0$ (see (\ref{eqq002}) in Section 4). In particular, if $\mu_1=\mu_2$,
then $\bb_0=\frac{2}{N-2}\max\{\mu_1, \mu_2\}$.

(2) For the case $p=2$ and $N\le 3$, some uniqueness results about system (\ref{eq2}) were introduced in \cite{WY}. However, their proofs heavily
depend on $p=2$, and so can not work here.\er

Now, let us consider the general case $-\la_1(\om)<\la_1, \la_2<0$.
Without loss of generality, we may assume that $\la_1\le \la_2$.

\bt\label{th1} Assume that $-\la_1(\om)<\la_1\le\la_2<0$.
Then system (\ref{eq2}) has a positive least energy solution $(u, v)$ with $E(u, v)=B$ for any $\beta\neq 0$.
\et

\br\label{remark} For the general case $\la_1\le \la_2$, when $N=4$ and $2p=2^\ast$, \cite[Theorem 1.3]{CZ} said that (\ref{eq2}) has a positive least energy
solution for any
$$\bb\subset (-\iy, 0)\cup(0, \bb_1)\cup(\bb_2, +\iy),$$
where $\bb_i, i=1, 2$ are some positive constants satisfying
$$\bb_1\le \min\{\mu_1,\mu_2\}\le \max\{\mu_1,\mu_2\}\le \bb_2.$$
That is, we do not know whether the least energy solution exists or not if $\bb\in [\bb_1, \bb_2]$ (In the symmetric case $\la_1=\la_2$, Remark \ref{rmk}-(1) says that nontrivial positive solutions do not exist if $\beta\in [\min\{\mu_1,\mu_2\}, \max\{\mu_1,\mu_2\}] $ and $\mu_1\neq \mu_2$).
Comparing this with Theorem \ref{th1}, it turns out that the general case $N\ge 5$ is completely different from
the special case $N=4$.
\er

Now, we study the limit behavior of the positive least energy solutions in the repulsive case $\beta\to -\infty$.
It is expected that
components of the limiting profile tend to separate in different regions of the underlying domain $\om$.
This phenomena, called {\it phase separation},
has been well studied for $L^{\iy}$ -bounded positive solutions of system (\ref{eq2}) in subcritical case $2p<2^\ast$ by \cite{WW2, WW3, NTTV}. The critical
case
$N=4$ and $p=2$ was studied by \cite{CZ}.
For other kinds of elliptic systems with strong competition, phase separation has also been well studied, we refer to \cite{CL0, CR, CTV} and references
therein.
Denote
$\{u>0\}:=\{x\in\om : u(x)>0 \}$ and $u^{\pm}:=\max\{\pm u, 0\}$. Then we have the following result.

\bt\label{th3} Assume that $-\la_1(\om)<\la_1\le\la_2<0$. Let $\beta_n<0,\,\,n\in\mathbb{N}$, satisfy $\beta_n\to-\iy$ as $n\to\iy$,
and $(u_n, v_n)$ be the positive least energy solutions of (\ref{eq2}) with $\beta=\beta_n$  which exists by Theorem \ref{th1}.
Then $\int_{\om}\beta_n u_n^p v_n^p\,dx\to 0$
as $n\to\iy$, and passing to a subsequence, one of the following conclusions holds.
\begin{itemize}

\item [(1)] $u_n\to u_{\iy}$ strongly in $H_0^1(\om)$ and $v_n\rightharpoonup 0$ weakly in $H_0^1(\om)$ (so $v_n\to 0$ for almost every $x\in \om$),
where $u_\iy$ is a positive least energy
solution of
$$-\Delta u +\la_1 u=\mu_1 |u|^{2^\ast-2}u,\,\, u\in H^1_0(\om).$$
\item [(2)] $v_n\to v_{\iy}$ strongly in $H_0^1(\om)$ and $u_n\rightharpoonup 0$ weakly in $H_0^1(\om)$ (so $u_n\to 0$ for almost every $x\in \om$),
where $v_\iy$ is a positive least energy
solution of
$$-\Delta v +\la_2 v=\mu_2 |v|^{2^\ast-2}v,\,\, v\in H^1_0(\om).$$

\item [(3)] $(u_n, v_n)\to (u_{\iy}, v_{\iy})$ strongly in $H$ and $u_\iy\cdot v_\iy\equiv 0$,
where $u_\iy\in C(\overline{\om})$ is a positive least energy
solution of
$$-\Delta u +\la_1 u=\mu_1 |u|^{2^\ast-2}u,\,\, u\in H^1_0(\{u_\iy>0\}),$$
and
$v_\iy\in C(\overline{\om})$ is a positive least energy
solution of
$$-\Delta v +\la_2 v=\mu_2 |v|^{2^\ast-2}v,\,\, v\in H^1_0(\{v_\iy>0\}).$$
Furthermore, both $\{v_\iy>0\}$ and $\{u_\iy>0\}$ are connected domains, and $\{v_\iy>0\}=\om\backslash\overline{\{u_\iy>0\}}$.

\end{itemize}
In particular, if $N\ge 6$, then only conclusion (3) holds, and $u_\iy-v_\iy$ is a least energy sign-changing solution to problem
\be\label{fc3}-\Delta u+\la_1 u^+-\la_2 u^-=\mu_1 (u^+)^{2^\ast-1}-\mu_2 (u^-)^{2^\ast-1}, \quad u\in H_0^1(\om).\ee\et

Here a sign-changing solution $u$ of (\ref{fc3}) is called {\it a least energy sign-changing solution}, if $u$ attains the minimal functional energy among all sign-changing solutions of (\ref{fc3}). As an application of Theorem \ref{th3}, we turn to consider the Brezis-Nirenberg problem
\be\label{equation1}-\Delta u +\la_1 u= \mu_1|u|^{2^\ast-2}u,\,\, u\in H^1_0(\om),\ee
where $-\la_1(\om)<\la_1<0$. Its corresponding functional is
$$J(u):=\frac{1}{2}\int_{\om}(|\nabla u|^2+\la_1 u^2)\,dx-\frac{1}{2^\ast}\int_{\om}\mu_1|u|^{2^\ast}\,dx,\,\,u\in H^1_0(\om).$$
Then we have the following result.

\bt\label{sign-changing} Assume $N\ge 6$. Let $(u_\iy, v_\iy)$ be in Theorem \ref{th3} in the symmetric
case where $\la_2=\la_1$ and $\mu_2=\mu_1$. Then $u_\iy-v_\iy$ is a
least energy sign-changing solution of (\ref{equation1}), and
\be\label{signchange} J(u_\iy-v_\iy)<B_{\mu_1}+\frac{1}{N}\mu_1^{-\frac{N-2}{2}}S^{\frac{N}{2}},\ee
where $B_{\mu_1}$ is the least energy of problem (\ref{equation1}) (see (\ref{eq11}) in Section 3).\et

\br\label{remark3}
\begin{itemize}

\item[$(i)$]
Theorem \ref{th3} has been proved for the special case $N=4$ and $2p=2^\ast$ by the authors(\cite{CZ}), where we raised an open question: Can one show that only conclusion (3) holds? The reviewer of \cite{CZ} pointed out that this question may be related to the existence of sign-changing solutions to the Brezis-Nirenberg problem (\ref{equation1}). Motivated by the reviewer's comment, it is natural for us to consider (\ref{eq2}) under assumption (\ref{fc}) in this paper.
Here in the case $N\ge 6$, we exclude conclusions (1)-(2) and verifies the reviewer's comment successfully, and so system (\ref{eq2}) is closely related to the Brezis-Nirenberg problem. Unfortunately, we do not know whether only Theorem \ref{th3}-(3) holds for $N=4,
5$, which still remains as an interesting open question.

\item[$(ii)$] In the proof of Theorem \ref{th3}-(3), a key point is to prove the continuity of $u_\iy$ and $v_\iy$. We remark here that, our proof of the continuity of $u_\iy$ and $v_\iy$ is completely different from that in \cite{CZ} for
the special case $N=4$, and can also be used to the special case $N=4$.

\item[$(iii)$] The existence of least energy sign-changing solutions to the Brezis-Nirenberg problem (\ref{equation1}) in the case $N\ge 6$ was proved in \cite{CSS} in 1986. Here, Theorem \ref{sign-changing} is a direct corollary of Theorem \ref{th3}, and so the proof of Theorem \ref{sign-changing} is completely different from \cite{CSS}.
\end{itemize}
\er

Since the nonlinearity and the coupling term are both critical in (\ref{eq2}),
the existence of nontrivial solutions of (\ref{eq2}) depends heavily
on the existence of the least energy solution of the following limit problem
{\allowdisplaybreaks
\be\label{eq4}
\begin{cases}-\Delta u  =
\mu_1 |u|^{2p-2}u+\beta |u|^{p-2}u|v|^p, & x\in \R^N,\\
-\Delta v  =\mu_2 |v|^{2p-2}v+\beta |v|^{p-2}v|u|^{p},    & x\in
\R^N,\\
u, v\in D^{1, 2}(\R^N),\end{cases}\ee
}%
where $D^{1,2}(\R^N):=\{u\in L^2(\R^N)\,:\, |\nabla u|\in L^2(\R^N)\}$ with
norm $\|u\|_{D^{1,2}}:=(\int_{\R^N}{|\nabla u|^2}\,dx)^{1/2}$. Let $S$ be the
sharp constant of $D^{1,2}(\R^N)\hookrightarrow L^{2^\ast}(\RN)$
\be\label{eq5}\intR{|\nabla u|^2}\ge
S\left(\intR{|u|^{2^\ast}}\right)^{\frac{2}{2^\ast}}.\ee

For $\varepsilon > 0$ and $y \in \RN$, we consider the Aubin-Talenti instanton
\cite{A, T} $U_{\varepsilon,y} \in D^{1,2}(\RN)$ defined by
\be\label{A-T} U_{\varepsilon,y}(x):= [N(N -
2)]^{\frac{N-2}{4}}\left(\frac{\varepsilon}{\varepsilon^2 +
|x-y|^2}\right)^{\frac{N-2}{2}}.\ee
Then $U_{\varepsilon,y}$ satisfies $-\Delta u = |u|^{2^\ast-2}u$ in $\RN$ and
\be\label{A-T1}\int_{\RN}|\nabla U_{\varepsilon,y}|^2\,dx =
\int_{\RN} |U_{\varepsilon,y}|^{2^\ast}\,dx = S^{N/2}.\ee Furthermore,
$\{U_{\varepsilon,y} : \varepsilon>0, y \in \R^N\}$ contains all positive solutions
of the equation $-\Delta u = |u|^{2^\ast-2}u$ in $\RN$.

Clearly (\ref{eq4}) has semi-trivial solutions $(\mu_1^{-\frac{N-2}{4}}U_{\e,y}, 0)$ and $(0, \mu_2^\frac{N-2}{4}U_{\e,y})$. Here, we are only
interested in nontrivial solutions of (\ref{eq4}). Define $D:=D^{1,2}(\RN)\times D^{1,2}(\RN)$
and a $C^1$ functional $I: D\to \R$
given by
\begin{align}\label{eq6}
I(u, v):=&\frac{1}{2}\int_{\RN}(|\nabla u|^2+|\nabla v|^2)-\frac{1}{2p}\int_{\RN}(\mu_1 |u|^{2p}+2\beta |u|^p|v|^p +\mu_2 |v|^{2p}).
\end{align}
As in \cite{LW1}, we consider the set
\begin{align*}\mathcal{N}=\left\{(u, v)\in D\,\,:\,\, u\not\equiv 0, \,v\not\equiv 0,\, I'(u, v)(u, 0)=I'(u, v)(0, v)=0\right\}.\end{align*}
Then any nontrivial solutions of (\ref{eq4}) belong to $\mathcal{N}$.
Similarly $\mathcal{N}\neq \emptyset$.
 We set
\be\label{eq7}A:=\inf_{(u, v)\in \mathcal{N}} I(u, v)=\inf_{(u, v)\in \mathcal{N}}\frac{1}{N}\intR{(|\nabla u|^2+|\nabla v|^2)}.\ee
Then we have the following theorem, which plays an important role in the proof of Theorem \ref{th1}.
\bt\label{th2}
\begin{itemize}

\item [(1)] If $\beta<0$, then $A$ is not attained.

\item [(2)] If $\bb>0$, then (\ref{eq4}) has a positive least energy solution $(U, V)$ with $I(U, V)=A$, which is radially symmetric decreasing. Moreover,

\begin{itemize}

\item [(2-1)]
if $\bb\ge \frac{2}{N-2}\max\{\mu_1, \mu_2\}$, then $I(\sqrt{k_0}U_{\e, y}, \sqrt{l_0}U_{\e, y})
=A$, where $(k_0, l_0)$ in (\ref{eqqqq}). That is, $(\sqrt{k_0}U_{\e, y}, \sqrt{l_0}U_{\e, y})$ is a positive least energy solution of (\ref{eq4}).

\item [(2-2)]
there exists $0<\bb_1\le\frac{2}{N-2}\max\{\mu_1, \mu_2\}$, and for any $0<\bb<\bb_1$, there exists a solution $(k(\bb), l(\bb))$ of (\ref{eq8}), such
that $$I(\sqrt{k(\bb)}U_{\e, y}, \sqrt{l(\bb)}U_{\e, y})
>A=I(U, V).$$ That is, $(\sqrt{k(\bb)}U_{\e, y}, \sqrt{l(\bb)}U_{\e, y})$ is a different positive solution of (\ref{eq4}) with respect to $(U, V)$.

\end{itemize}

\end{itemize}\et

\br\label{rmk1}In the case $N=4$ and $2p=2^\ast$, \cite[Theorem 1.5]{CZ} said that (\ref{eq4}) has no nontrivial nonnegative solution if
$\beta\in [\min\{\mu_1,\mu_2\}, \max\{\mu_1,\mu_2\}] $ and $\mu_1\neq \mu_2$;
$(\sqrt{k(\bb)}U_{\e, y}, \sqrt{l(\bb)}U_{\e, y})$ is a positive least energy solution of (\ref{eq4}) if $0<\bb<\min\{\mu_1,\mu_2\}$. Hence the general case $N\ge 5$ is completely different from
the case $N=4$. As we will see in Section 2,
the idea of proving Theorem \ref{th2}-(2) in case $0<\bb<\frac{2}{N-2}\max\{\mu_1, \mu_2\}$, which also works for the case $\bb\ge \frac{2}{N-2}\max\{\mu_1,
\mu_2\}$, is completely different from that in case $N=4$ (\cite{CZ}).\er

We can also study the uniqueness of the positive least energy solutions of (\ref{eq4}) just as Theorem \ref{unique1}.

\bt\label{unique2} Let $\bb_0$ be in Theorem \ref{unique1} and assume that $\bb>\bb_0$. Let $(u, v)$ be any a positive least energy solution
of (\ref{eq4}). Then $(u, v)=(\sqrt{k_0}U_{\e, y}, \sqrt{l_0}U_{\e, y})$ for some $\e>0$ and $y\in\RN$.\et

The rest of this paper proves these theorems, and some ideas of the proofs are similar to those in \cite{CZ}.
However, as pointed out above, the general case $N\ge 5$ is quite different from $N=4$,
and some new ideas are needed. We give some notations here. Throughout this paper,
we denote the norm of $L^q(\om)$ by $|u|_q =
(\int_{\om}|u|^q\,dx)^{\frac{1}{q}}$, the norm of $H^1_0(\om)$ by $\|u\|=|\nabla u|_2$ and positive constants
(possibly different) by $C$. The paper is organized as follows. Theorems \ref{th0} and \ref{th2} are proved
in Section 2, and we will see that these proofs are more delicate than those in case $N=4$ (\cite{CZ}).
In Section 3, we use Nehari manifold approach and Ekeland variational principle to prove Theorem \ref{th1} for the case $\bb<0$,
and use mountain pass argument to Theorem \ref{th1} for the case $\bb>0$. In Section 4, we use an elementary approach to prove Theorems \ref{unique1} and
\ref{unique2}. Finally, we use energy estimate methods to prove
Theorems \ref{th3} and \ref{sign-changing} in Section 5, where some different ideas are needed.

After the completion of this paper (see arXiv:1209.2522v1 for the original version), we learned that (\ref{eq2}) has also been studied in \cite{SK}, where the author showed the existence of
least energy solutions for $\bb>\bb_0>0$, where $\bb_0$ is an unknown constant. For $|\bb|$ sufficiently small, the existence of nontrivial solutions was studied in \cite{SK} via a perturbation method,
but the solutions obtained there seem not necessary to be least energy solutions.
Theorem \ref{th1} in this
paper is much more general that those results in \cite{SK}.
We also remark that any other results in our paper can not be found in \cite{SK}.

\vskip0.1in

\s{Proof of Theorems \ref{th0} and \ref{th2}}
\renewcommand{\theequation}{2.\arabic{equation}}

Define functions
{\allowdisplaybreaks
\begin{align}
\label{eqq01}&\al_1(k, l):=\mu_1 k^{p-1}+\bb k^{\frac{p}{2}-1}l^{\frac{p}{2}}-1, \quad k>0,\,\,l\ge 0;\\
\label{eqq02}&\al_2(k, l):=\mu_2 l^{p-1}+\bb l^{\frac{p}{2}-1}k^{\frac{p}{2}}-1,\quad l>0,\,\,k\ge 0;\\
\label{eqq03}&h_1(k):=\bb^{-2/p}\left(k^{1-p/2}-\mu_1 k^{p/2}\right)^{2/p}, \quad 0<k\le \mu_1^{-\frac{1}{p-1}};\\
\label{eqq04}&h_2(l):=\bb^{-2/p}\left(l^{1-p/2}-\mu_2 l^{p/2}\right)^{2/p}, \quad 0<l\le \mu_2^{-\frac{1}{p-1}}.
\end{align}
}%
Then $\al_1(k, h_1(k))\equiv 0$ and $\al_2(h_2(l), l)\equiv 0$.

\bl\label{lem}Assume that $\bb>0$, then equation
\be\label{eqq05}\al_1(k, l)=0,\quad \al_2(k, l)=0,\quad k, l>0\ee has a solution $(k_0, l_0)$, which satisfies
\be\label{eqq06}\al_2(k, h_1(k))<0,\quad \forall\,\,0<k<k_0,\ee
that is, $(k_0, l_0)$ satisfies (\ref{eqqqq}). Similarly, (\ref{eqq05}) has a solution $(k_1, l_1)$ such that
\be\label{eqq07}\al_1(h_2(l), l)<0,\quad \forall\,\,0<l<l_1.\ee
\el

\noindent {\bf Proof. } Equation $\al_1 (k, l)=0,\,\,k, l>0$ imply that
$$l=h_1(k), \quad 0<k<\mu_1^{-\frac{1}{p-1}}.$$
While, $\al_2(k, l)=0$  implies that $\mu_2 l^{p/2}+\beta  k^{p/2}  =l^{1-p/2}.$ Therefore, we turn to prove that
\be\label{eqq1}\mu_2\frac{1-\mu_1 k^{p-1}}{\beta  k^{p/2-1}}+\beta  k^{p/2}  =\left(\frac{1-\mu_1 k^{p-1}}{\beta  k^{p/2-1}}\right)^\frac{2-p}{p},\quad
0<k^{p-1}<\frac{1}{\mu_1}\ee
have a solution. Note that (\ref{eqq1}) is equivalent to
\be\label{eqq2}f(k):=\left(\frac{1}{\bb k^{p-1}}-\frac{\mu_1}{\bb}\right)^\frac{2-p}{p}-\frac{\mu_2}{\bb}-\frac{\bb^2-\mu_1\mu_2}{\bb}k^{p-1}=0,\,\,
0<k^{p-1}<\frac{1}{\mu_1}.\ee
Recall that $N\ge 5$ and $2p=2^\ast$, we have $2-p>0$ and so
$$\lim_{k\to 0+}f(k)=+\iy,\quad f(\mu_1^{-\frac{1}{p-1}})=-\frac{\bb}{\mu_1}<0.$$
Therefore, there exists $k_0\in (0, \mu_1^{-\frac{1}{p-1}})$ such that $f(k_0)=0$ and $f(k)>0$ for $k\in (0, k_0)$. Let $l_0=h_1(k_0)$, then
$(k_0, l_0)$ is a solution of (\ref{eqq05}). Moreover, (\ref{eqq06}) follows directly from $f(k)>0$ for $k\in (0, k_0)$. The existence of $(k_1, l_1)$ that
satisfy (\ref{eqq05}) and
(\ref{eqq07}) is similar.\hfill$\square$

\bl\label{lem3}Assume that $\bb\ge (p-1)\max\{\mu_1, \mu_2\}$, then $h_1(k)+k$ is strictly increasing for $k\in [0, \mu_1^{-\frac{1}{p-1}}]$ and
$h_2(l)+l$ is strictly increasing for $ l\in [0,\mu_2^{-\frac{1}{p-1}}].$ \el

\noindent {\bf Proof. } Since for $k>0$
$$h_1'(k)=\frac{2}{p}\bb^{-2/p}\left(k^{1-p/2}-\mu_1 k^{p/2}\right)^{2/p-1}\left((1-p/2)k^{-p/2}-\frac{p}{2}\mu_1 k^{p/2-1}\right),$$
we see that $h_1'(k)\ge 0$ for $0< \mu_1 k^{p-1}\le \frac{2-p}{p}$ or $\mu_1 k^{p-1}=1$, and $h_1'(k)<0$ for $\frac{2-p}{p}<\mu_1 k^{p-1}<1$.
By direct computations, we deduce from $h_1''(k)=0,\,\,\frac{2-p}{p}<\mu_1 k^{p-1}<1$
that $k=(\mu_1 p)^{-\frac{1}{p-1}}$. Since $\bb\ge (p-1)\max\{\mu_1, \mu_2\}$, we have
\begin{align*}
\min_{0<k^{p-1}\le\mu_1^{-1}}h_1'(k)=h_1'\Big((\mu_1 p)^{-\frac{1}{p-1}}\Big)=-\bb^{-2/p}\mu_1^{2/p}(p-1)^{2/p}\ge-1,
\end{align*}
and so $h_1'(k)>-1$ for $0<k\le\mu_1^{-\frac{1}{p-1}}$ with $k\neq (\mu_1 p)^{-\frac{1}{p-1}}$.
This implies that $h_1(k)+k$ is strictly increasing for $k\in [0, \mu_1^{-\frac{1}{p-1}}]$. Similarly,
$h_2(l)+l$ is strictly increasing for $ l\in [0,\mu_2^{-\frac{1}{p-1}}].$\hfill$\square$

\bl\label{lem4}Assume that $\bb\ge (p-1)\max\{\mu_1, \mu_2\}$. Let $(k_0, l_0)$ be in Lemma \ref{lem}. Then
$\max\{\mu_1(k_0+l_0)^{p-1}, \mu_2(k_0+l_0)^{p-1}\}<1$ and
\be\label{ff}\al_2(k, h_1(k))<0,\,\,\, \forall\,\,0<k<k_0;\quad\al_1(h_2(l), l)<0,\,\,\, \forall\,\,0<l<l_0.\ee\el

\noindent {\bf Proof. } By Lemma \ref{lem3} we have
$$h_1(\mu_1^{-\frac{1}{p-1}})+\mu_1^{-\frac{1}{p-1}}=\mu_1^{-\frac{1}{p-1}}>h_1(k_0)+k_0=k_0+l_0,$$
that is, $\mu_1(k_0+l_0)^{p-1}<1$. Similarly, $\mu_2(k_0+l_0)^{p-1}<1$. By Lemma \ref{lem}, to prove (\ref{ff}), it suffices to prove that $(k_0, l_0)=(k_1,
l_1)$.
By (\ref{eqq06})-(\ref{eqq07}) we see that $k_1\ge k_0, l_0\ge l_1$. If $k_1>k_0$, then $k_1+h_1(k_1)>k_0+h_1(k_0)$, that is,
$h_2(l_1)+l_1=k_1+l_1>k_0+l_0=h_2(l_0)+l_0$, and so $l_1>l_0$, a contradiction. Therefore, $k_1=k_0$ and $l_0=l_1$. This completes the proof.
\hfill$\square$

\bl\label{lem5}Assume that $\bb\ge (p-1)\max\{\mu_1, \mu_2\}$. Then
\be\label{eqq08}
\begin{cases}k+l \le k_0+l_0,
\\
\al_1(k, l)\ge 0,\quad\al_2(k, l)\ge 0,\\
k, l\ge 0,\,\,(k, l)\neq (0, 0)\end{cases}\ee
has a unique solution $(k_0, l_0)$.\el

\noindent {\bf Proof. } Note that $(k_0, l_0)$ satisfies (\ref{eqq08}). Let $(\tilde{k},\tilde{l})$ be any a solution of (\ref{eqq08}). Without loss of
generality, we assume that $\tilde{k}>0$. If $\tilde{l}=0$, then by $\tilde{k}\le k_0+l_0$ and $\al_1(\tilde{k}, 0)\ge 0$ we get that
$$1\le\mu_1 \tilde{k}^{p-1}\le \mu_1(k_0+l_0)^{p-1},$$
which contradicts with Lemma \ref{lem4}. Therefore $\tilde{l}>0$.

Assume by contradiction that $\tilde{k}< k_0$. Similarly to the proof of Lemma \ref{lem3},
by (\ref{eqq04}) it is easy to see that $h_2(l)$ is strictly increasing for $0<\mu_2 l^{p-1}\le \frac{2-p}{p}$, and strictly decreasing for
$\frac{2-p}{p}\le\mu_2 l^{p-1}\le 1$. Moreover, $h_2(0)=h_2(\mu_2^{-\frac{1}{p-1}})=0$. Since $0<\tilde{k}<k_0=h_2(l_0)$, there exists
$0<l_2<l_3<\mu_2^{-\frac{1}{p-1}}$ such that $h_2(l_2)=h_2(l_3)=\tilde{k}$ and
\be\label{eqq09} \al_2(\tilde{k}, l)<0\,\,\,\Longleftrightarrow\,\,\, h_2(l)>\tilde{k}\,\,\,\Longleftrightarrow\,\,\, l_2<l<l_3.\ee
Since $\al_2(\tilde{k}, \tilde{l})\ge 0$, we have $\tilde{l}\le l_2$ or $\tilde{l}\ge l_3$.
Since $\al_1(\tilde{k},\tilde{l})\ge 0$, we have $\tilde{l}\ge h_1(\tilde{k})$. By Lemma \ref{lem4} we have
$\al_2(\tilde{k}, h_1(\tilde{k}))<0$, and so $l_2<h_1(\tilde{k})<l_3$. These imply
\be\label{eqq010}\tilde{l}\ge l_3.\ee
On the other hand, since $l_1:=k_0+l_0-\tilde{k}>l_0$, we have
$$h_2(l_1)+k_0+l_0-\tilde{k}=h_2(l_1)+l_1>h_2(l_0)+l_0=k_0+l_0,$$
that is, $h_2(l_1)>\tilde{k}$. By (\ref{eqq09}) we have $l_2<l_1<l_3$. By $\tilde{k}+\tilde{l}\le k_0+l_0$ we have
$$\tilde{l}\le l_1<l_3,$$
which contradicts with (\ref{eqq010}). Therefore, $\tilde{k}\ge k_0$. By a similar argument, we also have $\tilde{l}\ge l_0$. Therefore,
$(\tilde{k}, \tilde{l})=(k_0, l_0)$. This completes the proof.\hfill$\square$\\

\noindent {\bf Proof of Theorem \ref{th0}. }
Assume that $-\la_1(\om)<\la_1=\la_2=\la<0$. By the Sobolev inequality (\ref{eq5}) it is standard to see that $B>0$.
Since
$\bb>0$, by Lemma \ref{lem} equation (\ref{eq8}) has a solution
$(k_0, l_0)$. Recall (\ref{func}), we see that $(\sqrt{k_0}\omega, \sqrt{l_0}\omega)$ is a nontrivial solution of (\ref{eq2}) and
\be\label{func2} 0<B\le E(\sqrt{k_0}\omega,\sqrt{l_0}\omega)=(k_0+l_0)B_1.\ee

Now we assume that $\bb\ge (p-1)\max\{\mu_1, \mu_2\}$, and we shall prove that $B=E(\sqrt{k_0}\omega,\sqrt{l_0}\omega)$. Let $\{(u_n, v_n)\}\subset\mathcal{M}$
be a minimizing sequence for $B$, that is, $E(u_n, v_n)\to B$. Define
$$c_n=\left(\int_{\om}|u_n|^{2p}\,dx\right)^{1/p},\quad d_n=\left(\int_{\om}|v_n|^{2p}\,dx\right)^{1/p}.$$
By (\ref{func1}) we have
{\allowdisplaybreaks
\begin{align}\label{eqq011}
   (N B_1)^{2/N}c_n\le\int_{\om}(|\nabla u_n|^2+\la u_n^2)&=\int_{\om}(\mu_1 |u_n|^{2p}+\beta |u_n|^p |v_n|^p)\nonumber\\
   &\le \mu_1 c_n^p+\beta c_n^{p/2} d_n^{p/2},\\
   \label{eqq012}(N B_1)^{2/N}d_n\le\int_{\om}(|\nabla v_n|^2+\la v_n^2)&=\int_{\om}(\mu_2 |v_n|^{2p}+\beta |u_n|^p |v_n|^p)\nonumber\\
   &\le \mu_2 d_n^p+\beta c_n^{p/2} d_n^{p/2}.
\end{align}
}%
Since $E(u_n, v_n)=\frac{1}{N}\int_{\om}(|\nabla u_n|^2+\la u_n^2+|\nabla v_n|^2+\la v_n^2)$,
by (\ref{func2}) we have
{\allowdisplaybreaks
\begin{gather}
    \label{eqq013}(N B_1)^{2/N}(c_n+d_n)\le N E(u_n,v_n)\le N(k_0+l_0)B_1+o(1),\\
    \label{eqq014}\mu_1 c_n^{p-1}+\beta c_n^{p/2-1} d_n^{p/2}\ge (N B_1)^{2/N},\\
    \label{eqq015}\mu_2 d_n^{p-1}+\beta c_n^{p/2} d_n^{p/2-1}\ge (N B_1)^{2/N}.
\end{gather}
}%
First, this means $c_n, d_n$ are uniformly bounded. Passing to a subsequence, we assume that $c_n\to c$ and $d_n\to d$.
Then by (\ref{eqq011})-(\ref{eqq012}) we
have $\mu_1 c^p+2\bb c^{p/2}d^{p/2}+\mu_2 d^p\ge NB>0$. Hence, without loss of generality, we assume that $c>0$. If $d=0$, then (\ref{eqq013}) implies
$c\le (NB_1)^{1-2/N}(k_0+l_0)$. By (\ref{eqq014}) and Lemma \ref{lem4} we get
$$(NB_1)^{2/N}\le\mu_1 c^{p-1}\le \mu_1 (k_0+l_0)^{p-1}(NB_1)^{2/N}<(NB_1)^{2/N},$$
a contradiction. Therefore, $c>0$ and $d>0$.
Let $k=\frac{c}{(N B_1)^{1-2/N}}$ and $l=\frac{d}{(N B_1)^{1-2/N}}$, then by (\ref{eqq013})-(\ref{eqq015}) we see that $(k, l)$ satisfies
(\ref{eqq08}). By Lemma \ref{lem5}
we see that $(k, l)=(k_0, l_0)$.
It follows that $c_n\to k_0 (N B_1)^{1-2/N}$ and $d_n\to  l_0 (N B_1)^{1-2/N}$ as $n\to +\iy$, and
$$NB=\lim_{n\to +\iy}NE(u_n, v_n)\ge \lim_{n\to +\iy}(N B_1)^{2/N}(c_n+d_n)=N(k_0+l_0)B_1.$$
Combining this with (\ref{func2}), one has that
$$ B =(k_0+l_0)B_1=E(\sqrt{k_0}\omega,\sqrt{l_0}\omega),$$
and so $(\sqrt{k_0}\omega,\sqrt{l_0}\omega)$ is a positive least energy solution of (\ref{eq2}).\hfill$\square$\\

Now we turn to the proof of Theorem \ref{th2}.
By the Sobolev inequality (\ref{eq5}) it is standard to see that
\be\label{equation}A=\inf_{(u, v)\in \mathcal{N}}\frac{1}{4}\intR{\big(|\nabla u|^2+|\nabla v|^2\big)}>0.\ee

\bl\label{lemma} If $A$ (resp. $B$) is attained by a couple $(u, v)\in \mathcal{N}$ (resp. $(u,v)\in\mathcal{M}$),
then this couple is a critical point of $I$ (resp. $E$),
provided $-\iy<\beta<0$. \el

\noindent {\bf Proof. } Let $\bb<0$. Assume that $(u, v)\in\mathcal{N}$ such that $A=I(u, v)$. Define
{\allowdisplaybreaks
\begin{gather*}
    G_1(u, v):=I'(u, v)(u, 0)=\int_{\RN}|\nabla u|^2
-\int_{\RN}(\mu_1 |u|^{2p}+\beta |u|^{p}|v|^p),\\
G_2(u, v):=I'(u, v)(0, v)=\int_{\RN}|\nabla v|^2-\int_{\RN}(\mu_2 |v|^{2p}+\beta |u|^p|v|^p).
\end{gather*}
}%
Then there exists two Lagrange multipliers $L_1, L_2\in\R$ such that
$$I'(u, v)+L_1 G_1'(u, v)+L_2G_2'(u, v)=0.$$
Acting on $(u, 0)$ and $(0, v)$ respectively, we obtain
\begin{align*}\left((2p-2)\int_{\RN}\mu_1 |u|^{2p}-(2-p)\int_{\RN}\bb |u|^p|v|^p\right)L_1+L_2p\int_{\RN}\bb |u|^p|v|^p=0,\\
\left((2p-2)\int_{\RN}\mu_2 |v|^{2p}-(2-p)\int_{\RN}\bb |u|^p|v|^p\right)L_2+L_1p\int_{\RN}\bb |u|^p|v|^p=0.
\end{align*}
Since $\bb<0$, we deduce from $G_1(u, v)=G_2(u, v)=0$ that
{\allowdisplaybreaks
\begin{align*}
&\left((2p-2)\int_{\RN}\mu_1 |u|^{2p}-(2-p)\int_{\RN}\bb |u|^p|v|^p\right)\\
&\quad\times\left((2p-2)\int_{\RN}\mu_2 |v|^{2p}-(2-p)\int_{\RN}\bb |u|^p|v|^p\right)>\left(p\int_{\RN}\bb |u|^p|v|^p\right)^2.
\end{align*}
}%
From this we have $L_1=L_2=0$ and so $I'(u, v)=0$. Similarly, if $(u, v)\in\mathcal{M}$ such that $E(u, v)=B$, then
$E'(u, v)=0$.
\hfill$\square$\\

\noindent {\bf Proof of (1) in Theorem \ref{th2}. } This proof is similar to the proof of \cite[Theorem 1.5-(1)]{CZ} in case $N=4$, but the details are
more delicate.
By (\ref{A-T}) we see that $\omega_{\mu_i}:=\mu_i^{-\frac{N-2}{4}}U_{1, 0}$ satisfies equation $-\Delta u=\mu_i |u|^{2^\ast-2}u$ in $\RN$.
Let $e_1=(1,0,\cdots,0)\in\RN$ and
$$(u_R(x), v_R(x))=(\omega_{\mu_1}(x), \omega_{\mu_2}(x+Re_1)).$$
Then $v_R\rightharpoonup 0$ weakly in $D^{1,2}(\RN)$ and so $v_R^p\rightharpoonup 0$ weakly in $L^{2}(\RN)$ as $R\to +\iy$.
That is,
$$\lim_{R\to+\iy}\intR{u^p_R v^p_R}=0.$$
Note that $\bb<0$. Then for $R>0$ sufficiently large, by a similar argument as that of Lemma \ref{lem} (or see the argument of existing $(t_\e, s_\e)$ in the
proof of Lemma \ref{lemma1} below), we see that
$$\begin{cases}t^2\intR{|\nabla u_R|^2}=t^2\mu_1\intR{u_R^{2p}}=t^{2p}\mu_1\intR{u_R^{2p}}+t^p s^p\beta\intR{u_R^pv_R^p},\\
s^2\intR{|\nabla v_R|^2}=s^2\mu_2\intR{v_R^{2p}}=s^{2p}\mu_2\intR{v_R^{2p}}+t^p s^p\beta\intR{u_R^pv_R^p},\end{cases}$$
have a solution $(t_{R}, s_{R})$ with $t_R>1$ and $s_R>1$. Denote
{\allowdisplaybreaks
\begin{gather*}
    D_1:=\mu_1\intR{u_R^{2p}}=\mu_1\intR{\omega_{\mu_1}^{2p}}>0,\\
    D_2:=\mu_2\intR{v_R^{2p}}=\mu_2\intR{\omega_{\mu_2}^{2p}}>0,\\
    F_R:=|\bb|\intR{u_R^pv_R^p}\to 0,\quad \hbox{as $R\to +\iy$}.
\end{gather*}
}%
Then
\be\label{fff}t_R^2 D_1=t_R^{2p}D_1-t_R^p s_R^p F_R,\quad s_R^2 D_2=s_R^{2p}D_2-t_R^p s_R^p F_R.\ee
Assume that, up to a subsequence, $t_R\to +\iy$ as $R\to \iy$, then by
$$t_R^{2p}D_1-t_R^2 D_1=s_R^{2p}D_2-s_R^2 D_2$$
we also have $s_R\to+\iy$. Note that $2-p<p$, we have
$$t_R^{p}D_1-t_R^{2-p} D_1\ge \frac{1}{2} t_R^{p}D_1,\,\,\,s_R^{p}D_2-s_R^{2-p} D_2\ge \frac{1}{2} s_R^{p}D_2,\,\,\hbox{for $R$ large enough,} $$
and so
$$F_R=\frac{t_R^{p}-t_R^{2-p} }{s_R^p}D_1\ge \frac{t_R^{p} }{2 s_R^p}D_1,\quad F_R=\frac{s_R^{p}-s_R^{2-p} }{t_R^p}D_2\ge \frac{s_R^{p} }{2 t_R^p}D_2,$$
which implies that
$$0<\frac{1}{4}D_1 D_2\le F_R^2\to 0, \quad \hbox{as $R\to +\iy$,}$$
a contradiction. Therefore, $t_R$ and $s_R$ are uniformly bounded. Then by (\ref{fff}) and $F_R\to 0$ as $R\to\iy$, we get that
$$\lim_{R\to+\iy}(|t_R-1|+|s_{R}-1|)=0.$$
Note that $(t_{R}u_R, s_{R} v_R)\in \mathcal{N}$, we see from (\ref{A-T1}) that
{\allowdisplaybreaks
\begin{align*}A&\le I(t_{R}u_R, s_R v_R)=\frac{1}{N}\left(t_R^2\intR{|\nabla u_R|^2}+s_R^2\intR{|\nabla v_R|^2}\right)\\
&=\frac{1}{N}\left(t_R^2\mu_1^{-\frac{N-2}{2}}+s_R^2\mu_2^{-\frac{N-2}{2}}\right)S^{N/2}.
\end{align*}
}%
Letting $R\to+\iy$, we get that $A\le\frac{1}{N}(\mu_1^{-\frac{N-2}{2}}+\mu_2^{-\frac{N-2}{2}})S^{N/2}$.

On the other hand, for any $(u, v)\in\mathcal{N}$, we see from $\beta<0$ and (\ref{eq5}) that
$$ \int_{\RN}|\nabla u|^2\,dx
\le \mu_1\int_{\RN} |u|^{2p}\,dx\le \mu_1 S^{-p}\left(\intR{|\nabla u|^2}\right)^{p},$$
and so $\int_{\RN}|\nabla u|^2\,dx\ge \mu_1^{-\frac{N-2}{2}}S^{N/2}$. Similarly, $\int_{\RN}|\nabla v|^2\,dx\ge \mu_2^{-\frac{N-2}{2}}S^{N/2}$.
Combining these with (\ref{eq7}), we get that $A\ge\frac{1}{N}(\mu_1^{-\frac{N-2}{2}}+\mu_2^{-\frac{N-2}{2}})S^{N/2}$. Hence,
\be\label{eq501}A=\frac{1}{N}\left(\mu_1^{-\frac{N-2}{2}}+\mu_2^{-\frac{N-2}{2}}\right)S^{N/2}.\ee

Now, assume that $A$ is attained by some $(u, v)\in \mathcal{N}$, then $(|u|, |v|)\in \mathcal{N}$ and $I(|u|, |v|)=A$.
By Lemma \ref{lemma}, we get that $(|u|, |v|)$ is a nontrivial solution of (\ref{eq4}). By the
maximum principle, we may assume that $u>0, v>0$ and so $\intR{u^p v^p}>0$. That is,
$$ \int_{\RN}|\nabla u|^2\,dx
< \mu_1\int_{\RN} |u|^{2p}\,dx\le \mu_1 S^{-p}\left(\intR{|\nabla u|^2}\right)^{p},$$
Therefore, it is easy to see that
$$A=I(u, v)=\frac{1}{N}\intR{(|\nabla u|^2+|\nabla v|^2)}>\frac{1}{N}\left(\mu_1^{-\frac{N-2}{2}}+\mu_2^{-\frac{N-2}{2}}\right)S^{N/2},$$
which is a contradiction. This completes the proof.\hfill$\square$\\

\noindent {\bf Proof of (2-1) in Theorem \ref{th2}. }This proof is similar to the proof of Theorem \ref{th0}.
Since $\bb>0$, by Lemma \ref{lem} equation (\ref{eq8}) has a solution
$(k_0, l_0)$. Then $(\sqrt{k_0}U_{\e,y},\sqrt{l_0}U_{\e,y})$ is a nontrivial solution of (\ref{eq4}) and
\be\label{ffc6}A\le I\left(\sqrt{k_0}U_{\e,y}, \sqrt{l_0}U_{\e,y}\right)=\frac{1}{N}(k_0+l_0)S^{N/2}.\ee
Assume that $\bb\ge (p-1)\max\{\mu_1, \mu_2\}$. Let $\{(u_n, v_n)\}\subset\mathcal{N}$ be a minimizing sequence for $A$, that is, $I(u_n, v_n)\to A$. Define
$c_n=\left(\intR{|u_n|^{2p}}\right)^{1/p}$,  $d_n=\left(\intR{|v_n|^{2p}}\right)^{1/p}$, we have
{\allowdisplaybreaks
\begin{gather*}
    Sc_n\le\intR{|\nabla u_n|^2}=\intR{\mu_1 |u_n|^{2p}+\beta |u_n|^p |v_n|^p}\le \mu_1 c_n^p+\beta c_n^{p/2} d_n^{p/2},\\
    Sd_n\le\intR{|\nabla v_n|^2}=\intR{\mu_2 |v_n|^{2p}+\beta |u_n|^p |v_n|^p}\le \mu_2 d_n^p+\beta c_n^{p/2} d_n^{p/2}.
\end{gather*}
}%
This means
{\allowdisplaybreaks
\begin{gather*}
    S(c_n+d_n)\le NI(u_n,v_n)\le (k_0+l_0)S^{N/2}+o(1),\\
    \mu_1 c_n^{p-1}+\beta c_n^{p/2-1} d_n^{p/2}\ge S,\quad
    \beta c_n^{p/2}d_n^{p/2-1} +\mu_2 d_n^{p-1}\ge S.
\end{gather*}
}%
Similarly as in the proof of Theorem \ref{th0}, we see that $c_n\to k_0 S^{N/2-1}$ and $d_n\to l_0 S^{N/2-1}$ as $n\to +\iy$, and
$$NA=\lim_{n\to +\iy}NI(u_n, v_n)\ge \lim_{n\to +\iy}S(c_n+d_n)=(k_0+l_0)S^{N/2}.$$
This implies that
\be\label{eq9}A=\frac{1}{N}(k_0+l_0)S^{N/2}=I(\sqrt{k_0}U_{\e,y},\sqrt{l_0}U_{\e,y}),\ee
and so $(\sqrt{k_0}U_{\e,y}, \sqrt{l_0}U_{\e,y})$ is a positive least energy solution of (\ref{eq4}).\hfill$\square$\\

To finish the proof of Theorem \ref{th2}, we need to show that (\ref{eq4}) has a positive least energy solution for any $0<\bb<(p-1)\max\{\mu_1, \mu_2\}$.
The following proof works for all $\bb>0$. Therefore, we assume that $\bb>0$, and define
\be\label{eqqq40} A':=\inf_{(u, v)\in\mathcal{N}'}I(u, v),\ee
where
\begin{align}\label{eqqq41}\mathcal{N}':=\left\{(u, v)\in D\setminus\{ (0,0)\},\,\,\, I'(u, v)(u, v)=0\right\}. \end{align}
Note that $\mathcal{N}\subset\mathcal{N}'$, one has that
$A'\le A$. By Sobolev inequality, we have $A'>0$. Define $B(0, R):=\{x\in\RN : |x|< R\}$ and $H(0, R):= H_0^1(B(0, R))\times H_0^1(B(0, R))$.
Consider

{\allowdisplaybreaks
\be\label{eqqq42}
\begin{cases}-\Delta u  =
\mu_1 |u|^{2p-2}u+\beta |u|^{p-2}u|v|^p, & x\in B(0,R),\\
-\Delta v  =\mu_2 |v|^{2p-2}v+\beta |v|^{p-2}v|u|^{p},    & x\in
B(0, R),\\
u, v\in H^1_0(B(0, R)),\end{cases}\ee
}%
and define
\be\label{eqqq43} A'(R):=\inf_{(u, v)\in\mathcal{N}'(R)}I(u, v),\ee
where
{\allowdisplaybreaks
\begin{align}\label{eqqq44}\mathcal{N}'(R) :=\Bigg\{&(u, v) \in H(0, R)\setminus\{ (0,0)\},\,\, \int_{B(0, R)}
(|\nabla u|^2+|\nabla v|^2)\nonumber\\
&
- \int_{B(0, R)}(\mu_1 |u|^{2p}+2\beta |u|^p|v|^p+\mu_2 |v|^{2p})=0\Bigg\}. \end{align}
}%

\bl\label{lemm1} $A'(R)\equiv A'$ for all $R>0$.\el

\noindent{\bf Proof. } Take any $R_1>R_2$. By $\mathcal{N}'(R_2)\subset \mathcal{N}'(R_1)$, we have $A'(R_1)\le A'(R_2)$. On the other hand, for any
$(u, v)\in \mathcal{N}'(R_1)$, we define $$(u_1(x), v_1(x)):=\left(\left(\frac{R_1}{R_2}\right)^{\frac{N-2}{2}}u\left(\frac{R_1}{R_2}x\right), \,\,
\left(\frac{R_1}{R_2}\right)^{\frac{N-2}{2}}v\left(\frac{R_1}{R_2}x\right)\right),$$
then it is standard to see that $(u_1, v_1)\in\mathcal{N}'(R_2)$, and so
$$A'(R_2)\le I(u_1, v_1)=I(u, v),\quad\forall (u, v)\in \mathcal{N}'(R_1).$$
That is, $A'(R_2)\le A'(R_1)$ and so $A'(R_1)=A'(R_2)$.

Clearly $A'\le A'(R)$. Let $(u_n, v_n)\in \mathcal{N}'$ be a minimizing sequence of $A'$. Moreover, we may assume that $u_n, v_n\in H_0^1(B(0, R_n))$
for some $R_n>0$. Then $(u_n, v_n)\in\mathcal{N}'(R_n)$ and
$$A'=\lim_{n\to\iy}I(u_n, v_n)\ge \lim_{n\to\iy}A'(R_n)\equiv A'(R).$$
Therefore, $A'(R)\equiv A'$ for all $R>0$.\hfill$\square$\\

Let $0\le \e<p-1$. Consider
{\allowdisplaybreaks
\be\label{eqqq45}
\begin{cases}-\Delta u  =
\mu_1 |u|^{2p-2-2\e}u+\beta |u|^{p-2-\e}u|v|^{p-\e}, & x\in B(0,1),\\
-\Delta v  =\mu_2 |v|^{2p-2-2\e}v+\beta |v|^{p-2-\e}v|u|^{p-\e},    & x\in
B(0, 1),\\
u, v\in H^1_0(B(0, 1)),\end{cases}\ee
}%
and define
\be\label{eqqq46} A_\e:=\inf_{(u, v)\in\mathcal{N}'_\e}I_\e(u, v),\ee
where
{\allowdisplaybreaks
\begin{align*}I_\e(u, v):=&\frac{1}{2}\int_{B(0,1)}(|\nabla u|^2+|\nabla v|^2)\\
&-\frac{1}{2p-2\e}\int_{B(0,1)}(\mu_1 |u|^{2p-2\e}+2\beta |u|^{p-\e}|v|^{p-\e} +\mu_2 |v|^{2p-2\e}),\\
\mathcal{N}'_\e :=&\left\{(u, v) \in H(0, 1)\setminus\{ (0,0)\},\,\, H_\e(u, v):=I_\e'(u, v)(u, v)=0\right\}.\end{align*}
}%

\bl\label{lemm2} For any $0<\e<p-1$, there holds
$$A_\e<\min\left\{\inf_{(u, 0)\in\mathcal{N}'_\e}I_\e(u, 0),\,\,\inf_{(0, v)\in\mathcal{N}'_\e}I_\e(0, v)\right\}.$$\el

\noindent{\bf Proof. } Fix any $0<\e<p-1$. Recall that $2<2p-2\e<2^\ast$, we may let $u_i$ be a least energy solution of
$$-\Delta u=\mu_i |u|^{2p-2-2\e}u,\quad u\in H_0^1(B(0, 1)),$$
Then $$I_\e(u_1, 0)=c_1:=\inf_{(u, 0)\in\mathcal{N}'_\e}I_\e(u, 0),\quad I_\e(0,u_2)=c_2:=\inf_{(0, v)\in\mathcal{N}'_\e}I_\e(0, v).$$
The following proof is inspired by \cite{AFP}. For any $s\in\R$, there exists a unique $t(s)>0$ such that $(t(s)u_1, t(s) su_2)\in \mathcal{N}'_\e$. In fact,
{\allowdisplaybreaks
\begin{align*}
t(s)^{2p-2\e-2}&=\frac{\int_{B(0,1)}(|\nabla u_1|^2+s^2|\nabla u_2|^2)}{\int_{B(0,1)}(\mu_1 |u_1|^{2p-2\e}+2\beta |u_1|^{p-\e}|su_2|^{p-\e}+\mu_2
|su_2|^{2p-2\e})}\\
&=\frac{p'c_1+s^2 p' c_2}{p'c_1+|s|^{2p-2\e}p'c_2+|s|^{p-\e}\int_{B(0, 1)}2\bb |u_1|^{p-\e}|u_2|^{p-\e}},
\end{align*}
}%
where $p'=\frac{2p-2\e}{p-1-\e}$. Note that $t(0)=1$. Recall that $1<p-\e<2$, by direct computations we have
\begin{align*}
\lim_{s\to 0}\frac{t'(s)}{|s|^{p-2-\e}s}=-\frac{(p-\e)\int_{B(0, 1)}2\bb |u_1|^{p-\e}|u_2|^{p-\e}}{(2p-2\e-2)p'c_1},
\end{align*}
that is,
$$t'(s)=-\frac{(p-\e)\int_{B(0, 1)}2\bb |u_1|^{p-\e}|u_2|^{p-\e}}{(2p-2\e-2)p'c_1}|s|^{p-\e-2}s(1+o(1)),\,\,\,\hbox{as $s\to 0$,}$$
and so
$$t(s)=1-\frac{\int_{B(0, 1)}2\bb |u_1|^{p-\e}|u_2|^{p-\e}}{(2p-2\e-2)p'c_1}|s|^{p-\e}(1+o(1)),\,\,\,\hbox{as $s\to 0$.}$$
This implies that
{\allowdisplaybreaks
\begin{align*}
t(s)^{2p-2\e}&=1-\frac{(2p-2\e)\int_{B(0, 1)}2\bb |u_1|^{p-\e}|u_2|^{p-\e}}{(2p-2\e-2)p'c_1}|s|^{p-\e}(1+o(1))\\
&=1-\frac{\int_{B(0, 1)}2\bb |u_1|^{p-\e}|u_2|^{p-\e}}{2c_1}|s|^{p-\e}(1+o(1)),\,\,\,\hbox{as $s\to 0$.}
\end{align*}
}%
Therefore, we see from $1/2-1/p'=1/(2p-2\e)>0$ that
{\allowdisplaybreaks
\begin{align*}
A_\e &\le I_\e\left(t(s)u_{1}, t(s)su_{2}\right)\\
&=\frac{t(s)^{2p-2\e}}{p'}\left(p'c_1+|s|^{2p-2\e}p'c_2+|s|^{p-\e}\int_{B(0, 1)}2\bb |u_1|^{p-\e}|u_2|^{p-\e}\right)\\
&=c_1-\left(\frac{1}{2}-\frac{1}{p'}\right)|s|^{p-\e}\int_{B(0, 1)}2\bb |u_1|^{p-\e}|u_2|^{p-\e}+o(|s|^{p-\e})\\
&<c_1=\inf_{(u, 0)\in\mathcal{N}'_\e}I_\e(u, 0)\quad \hbox{as $|s|>0$ small enough,}
\end{align*}
}%
By a similar argument, we have $A_\e<\inf\limits_{(0, v)\in\mathcal{N}'_\e}I_\e(0, v)$.\hfill$\square$\\

Recalling $\omega_{\mu_i}$ in the proof of Theorem \ref{th2} -(1), similarly as Lemma \ref{lemm2}, we have
{\allowdisplaybreaks
\begin{align}\label{eqqq47}A'<&\min\left\{\inf_{(u, 0)\in\mathcal{N}'}I(u, 0),\,\,\inf_{(0, v)\in\mathcal{N}'}I(0, v)\right\}
=\min\Big\{I(\omega_{\mu_1}, 0),\,I(0, \omega_{\mu_2})\Big\}\nonumber\\
=&\min\left\{\frac{1}{N}\mu_2^{-\frac{N-2}{2}} S^{N/2},\,\,\frac{1}{N}\mu_1^{-\frac{N-2}{2}} S^{N/2}\right\}.\end{align}
}%

\bt\label{lemm4} For any $0<\e<p-1$, (\ref{eqqq45}) has a classical least energy solution $(u_\e, v_\e)$, and $u_\e, v_\e$ are both positive radially symmetric
decreasing.\et

\noindent{\bf Proof. } Fix any $0<\e<p-1$, it is easy to see that $A_\e>0$. For $(u, v)\in \mathcal{N}'_\e$ with $u\ge 0, v\ge 0$,  we denote by $(u^*, v^*)$ as
its Schwartz
symmetrization. Then by the properties of Schwartz
symmetrization and $\bb>0$, we have
$$\int_{B(0, 1)}
(|\nabla u^*|^2+|\nabla v^*|^2)
\le \int_{B(0, 1)}(\mu_1 |u^*|^{2p-2\e}+2\beta |u^*|^{p-\e}|v^*|^{p-\e}+\mu_2 |v^*|^{2p-2\e}).$$
Therefore, there exists $0<t^*\le 1$ such that $(t^*u^*, t^* v^*)\in\mathcal{N}'_\e$, and then
{\allowdisplaybreaks
\begin{align}\label{eqqq48}
I_\e(t^*u^*, t^* v^*)&=\left(\frac{1}{2}-\frac{1}{2p-2\e}\right)(t^*)^2\int_{B(0, 1)}
(|\nabla u^*|^2+|\nabla v^*|^2)\nonumber\\
&\le\left(\frac{1}{2}-\frac{1}{2p-2\e}\right)\int_{B(0, 1)}
(|\nabla u|^2+|\nabla v|^2)=I_\e(u, v).
\end{align}
}%
Therefore, we may take a minimizing sequence $(u_n, v_n)\in\mathcal{N}'_\e$ of $A_\e$ such that $(u_n, v_n)=(u_n^*, v_n^*)$ and $I_\e(u_n, v_n) \to A_\e$.
We see from (\ref{eqqq48}) that $u_n, v_n$ are uniformly bounded in $H_0^1(B(0,1))$. Passing to a subsequence, we may assume that
$u_n\rightharpoonup u_\e, v_n\rightharpoonup v_\e$ weakly in $H_0^1(B(0,1))$. By the compactness of the embedding $H_0^1(B(0,1))\hookrightarrow
L^{2p-2\e}(B(0,1))$, we have
{\allowdisplaybreaks
\begin{align*}
&\int_{B(0, 1)}(\mu_1 |u_\e|^{2p-2\e}+2\beta |u_\e|^{p-\e}|v_\e|^{p-\e}+\mu_2 |v_\e|^{2p-2\e})\\
=&\lim_{n\to\iy}\int_{B(0, 1)}(\mu_1 |u_n|^{2p-2\e}+2\beta |u_n|^{p-\e}|v_n|^{p-\e}+\mu_2 |v_n|^{2p-2\e})\\
=&\frac{2p-2\e}{p-1-\e}\lim_{n\to\iy}I_\e(u_n, v_n)=\frac{2p-2\e}{p-1-\e}A_\e>0,
\end{align*}
}%
which implies $(u_\e, v_\e)\neq (0, 0)$. Moreover, $u_\e\ge 0, v_\e\ge 0$ are radially symmetric.
Meanwhile, $\int_{B(0, 1)}
(|\nabla u_\e|^2+|\nabla v_\e|^2)\le\lim\limits_{n\to\iy}\int_{B(0, 1)}
(|\nabla u_n|^2+|\nabla v_n|^2)$, so
$$\int_{B(0, 1)}
(|\nabla u_\e|^2+|\nabla v_\e|^2)
\le \int_{B(0, 1)}(\mu_1 u_\e^{2p-2\e}+2\beta u_\e^{p-\e}v_\e^{p-\e}+\mu_2 v_\e^{2p-2\e}).$$
Therefore, there exists $0<t_\e\le 1$ such that $(t_\e u_\e, t_\e v_\e)\in\mathcal{N}'_\e$, and then
{\allowdisplaybreaks
\begin{align*}
A_\e&\le I_\e(t_\e u_\e, t_\e v_\e)=\left(\frac{1}{2}-\frac{1}{2p-2\e}\right)(t_\e)^2\int_{B(0, 1)}
(|\nabla u_\e|^2+|\nabla v_\e|^2)\\
&\le\lim_{n\to\iy}\left(\frac{1}{2}-\frac{1}{2p-2\e}\right)\int_{B(0, 1)}
(|\nabla u_n|^2+|\nabla v_n|^2)\\
&=\lim_{n\to\iy} I_\e(u_n, v_n)=A_\e.
\end{align*}
}%
Therefore, $t_\e=1$ and $(u_\e, v_\e)\in\mathcal{N}'_\e$ with $I(u_\e, v_\e)=A_\e$. Moreover,
$$\int_{B(0, 1)}
(|\nabla u_\e|^2+|\nabla v_\e|^2)=\lim\limits_{n\to\iy}\int_{B(0, 1)}
(|\nabla u_n|^2+|\nabla v_n|^2),$$
that is, $u_n\to u_\e$ and $v_n\to v_\e$ strongly in $H_0^1(B(0, 1))$.
There exists a Lagrange multiplier
$\gamma\in \R$ such that
$$I_\e'(u_\e, v_\e)-\gamma H'_\e(u_\e, v_\e)=0.$$
Since $I'_\e(u_\e, v_\e)(u_\e, v_\e)=H_\e(u_\e, v_\e)=0$ and $$H_\e'(u_\e, v_\e)(u_\e, v_\e)=(2+2\e-2p)\int_{B(0, 1)}(\mu_1 u_\e^{2p-2\e}+2\beta
u_\e^{p-\e}v_\e^{p-\e}+\mu_2 v_\e^{2p-2\e})< 0,$$
we get that $\gamma=0$ and so $I'_\e(u_\e, v_\e)=0$. By Lemma \ref{lemm2}, we see that $u_\e\not\equiv 0$ and $v_\e\not\equiv 0$.
This means that
$(u_\e, v_\e)$ is a least energy solution of (\ref{eqqq43}). Recall that $u_\e, v_\e\ge 0$ are radially symmetric non-increasing. By regularity theory and the
maximum principle, we see that $u_\e, v_\e>0$ in $B(0, 1)$, $u_\e, v_\e\in C^2(B(0,1))$ and are radially symmetric decreasing.
\hfill$\square$\\

\noindent {\bf Completion of the proof of (2) in Theorem \ref{th2}. }Recalling (\ref{eqqq44}), for any $(u, v)\in \mathcal{N}'(1)$, there exists
$t_\e>0$ such that $(t_\e u, t_\e v)\in \mathcal{N}'_\e$ with $t_\e\to 1$ as $\e\to 0$. Then
$$\limsup_{\e\to 0}A_\e\le\limsup_{\e\to 0}I_\e(t_\e u, t_\e v)=I(u, v), \,\,\,\forall (u, v)\in \mathcal{N}'(1).$$
By Lemma \ref{lemm1} we have \be\label{eqqq52}\limsup_{\e\to 0}A_\e\le A'(1)=A'.\ee
By Theorem \ref{lemm4}, let $(u_\e, v_\e)$ be a positive
least energy solution of (\ref{eqqq45}), which is radially symmetric decreasing. By
$I_\e'(u_\e, v_\e)(u_\e, v_\e)=0$ and Sobolev inequality, it is easily seen that
\be\label{eqqq50}\frac{2p-2\e}{p-\e-1}A_\e=\int_{B(0, 1)}(|\nabla u_\e|^2+|\nabla v_\e|^2)\ge C_0,\quad \forall\,\, 0<\e\le \frac{p-1}{2},\ee
where $C_0$ is a positive constant independent of $\e$. Then
$u_\e, v_\e$ are uniformly bounded in $H_0^1(B(0, 1))$. Passing to a subsequence, we may assume that
$u_\e\rightharpoonup u_0$ and $v_\e\rightharpoonup v_0$ weakly in $H_0^1(B(0,1))$. Then $(u_0, v_0)$ is a solution of
{\allowdisplaybreaks
\be\label{eqqq49}
\begin{cases}-\Delta u  =
\mu_1 |u|^{2p-2}u+\beta |u|^{p-2}u|v|^{p}, & x\in B(0,1),\\
-\Delta v  =\mu_2 |v|^{2p-2}v+\beta |v|^{p-2}v|u|^{p},    & x\in
B(0, 1),\\
u, v\in H^1_0(B(0, 1)).\end{cases}\ee
}%

Assume by contradiction that $\|u_\e\|_{\iy}+\|v_\e\|_{\iy}$ is uniformly bounded, then by the Dominated Convergent Theorem, we get that
{\allowdisplaybreaks
\begin{gather*}
    \lim_{\e\to 0}\int_{B(0, 1)}u_\e^{2p-2\e}=\int_{B(0, 1)}u_0^{2p},\quad \lim_{\e\to 0}\int_{B(0, 1)}v_\e^{2p-2\e}=\int_{B(0, 1)}v_0^{2p},\\
    \lim_{\e\to 0}\int_{B(0, 1)}u_\e^{p-\e}v_\e^{p-\e}=\int_{B(0, 1)}u_0^{p}v_0^p.
\end{gather*}
}%
Combining these with $I'_\e(u_\e, v_\e)=I'(u_0, v_0)=0$, it is standard to show that $u_\e\to u_0$ and $v_\e\to v_0$ strongly in $H_0^1(B(0, 1))$.
Then by
(\ref{eqqq50}), we see that $(u_0, v_0)\neq (0, 0)$. Moreover, $u_0\ge 0, v_0\ge 0$. We may assume that $u_0\not\equiv 0$. By the strong maximum principle,
$u_0>0$ in $B(0,1)$. Note that $2p=2^\ast$. Combining these with Pohozaev identity, we have
$$0<\int_{\partial B(0,1)}(|\nabla u_0|^2+ |\nabla v_0|^2)(x\cdot \nu)\,d\sigma= 0,$$
a contradiction. Here, $\nu$ denotes the outward unit normal vector on $\partial B(0, 1)$. Therefore, $\|u_\e\|_{\iy}+\|v_\e\|_{\iy}\to \iy$ as $\e\to 0$. We
will use a blowup analysis. Note that $u_\e(0)=\max\limits_{B(0,1)}u_\e(x)$ and $v_\e(0)=\max\limits_{B(0,1)}v_\e(x)$, we define $K_\e:=max\{u_\e (0),
v_\e(0)\}$, then $K_\e\to
+\iy$. Define
$$U_\e(x)=K_\e^{-1}u_\e(K_\e^{-\al_\e} x),\quad V_\e(x)=K_\e^{-1}v_\e(K_\e^{-\al_\e} x),\quad \al_\e=p-1-\e.$$
Then
\be\label{eqqq51}1=\max\{U_\e (0), V_\e(0)\}=\max\left\{\max\limits_{x\in B(0, K_\e^{\al_\e})}U_\e(x),
\max\limits_{x\in B(0, K_\e^{\al_\e})}V_\e(x)\right\}\ee and $U_\e, V_\e$ satisfy
\begin{displaymath}
\begin{cases}-\Delta U_\e =
\mu_1 U_\e^{2p-2\e-1}+\bb U_\e^{p-1-\e} V_\e^{p-\e}, & x\in B(0, K_\e^{\al_\e}),\\
-\Delta V_\e =
\mu_2 V_\e^{2p-2\e-1}+\bb V_\e^{p-1-\e} U_\e^{p-\e}, & x\in B(0, K_\e^{\al_\e}).\end{cases}
\end{displaymath}
Since
$$\intR{|\nabla U_\e|^2}=K_\e ^{-(N-2)\e}\intR{|\nabla
u_\e|^2}\le\intR{|\nabla
u_\e|^2},$$ we see that $\{(U_\e, V_\e)\}_{n\ge 1}$ is bounded in
$D^{1,2}(\RN)\times D^{1,2}(\RN)= D$. By elliptic estimates, for a
subsequence we have $(U_\e, V_\e)\to (U, V)\in D$ uniformly in every
compact subset of $\RN$ as $\e\to 0$, and $(U, V)$ satisfies (\ref{eq4}), that is $I'(U, V)=0$.
Moreover, $U\ge 0, V\ge 0$ are radially symmetric non-increasing. By (\ref{eqqq51}) we have $(U, V)\neq (0, 0)$, and so $(U, V)\in \mathcal{N}'$. Then we deduce
from (\ref{eqqq52})
that
{\allowdisplaybreaks
\begin{align*}
A' &\le I(U, V)=\left(\frac{1}{2}-\frac{1}{2p}\right)\int_{\RN}(|\nabla U|^2+|\nabla V|^2)\,dx\\
&\le\liminf_{\e\to 0}\left(\frac{1}{2}-\frac{1}{2p-2\e}\right)\int_{B(0, K_\e^{\al_\e})}(|\nabla U_\e|^2+|\nabla V_\e|^2)\,dx\\
&\le\liminf_{\e\to 0}\left(\frac{1}{2}-\frac{1}{2p-2\e}\right)\int_{B(0, 1)}(|\nabla u_\e|^2+|\nabla v_\e|^2)\,dx\\
&=\liminf_{\e\to 0}A_\e\le A'.
\end{align*}
}%
This implies that $I(U, V)=A'$. By (\ref{eqqq47}) we have that $U\not\equiv 0$ and $V\not\equiv 0$. By the strong maximum principle, $U>0$ and $V>0$ are radially
symmetric decreasing.
We also have $(U, V)\in \mathcal{N}$, and so $I(U, V)\ge A\ge A'$, that is,
\be\label{ffc2}I(U, V)=A=A',\ee and $(U, V)$ is a positive least energy solution of (\ref{eq4}), which is radially symmetric decreasing.

Finally, we show the existence of $(k(\bb), l(\bb))$ for $\bb>0$ small. Recall (\ref{eqq01})-(\ref{eqq02}),
We denote $\al_i(k, l)$ by $\al_i(k, l, \bb)$ here.
Define $k(0)=\mu_1^{-\frac{1}{p-1}}$ and $l(0)=\mu_2^{-\frac{1}{p-1}}$, then $\al_i(k(0), l(0), 0)=0, i=1, 2$.  Note that
{\allowdisplaybreaks
\begin{align*}
    \partial_k \al_1(k(0), l(0), 0)&=(p-1)\mu_1 k(0)^{p-2}>0,\\
   \partial_l \al_2(k(0), l(0), 0)&=(p-1)\mu_2 l(0)^{p-2}>0,\\
   \partial_l \al_1(k(0), l(0), 0)&=\partial_k \al_2(k(0), l(0), 0)=0,
\end{align*}
}%
which implies that
$$\det\left(
                     \begin{array}{cc}
                       \partial_k \al_1(k(0), l(0), 0) & \partial_l \al_1(k(0), l(0), 0) \\
                       \partial_k \al_2(k(0), l(0), 0) & \partial_l \al_2(k(0), l(0), 0) \\
                     \end{array}
                   \right)>0.
$$
Therefore, by the implicit function theorem, $k(\bb), l(\bb)$ are well defined and class $C^1$ on $(-\bb_2, \bb_2)$ for some $\bb_2>0$, and
$\al_i(k(\bb), l(\bb), \bb)\equiv 0, i=1,2$. This implies that $(\sqrt{k(\bb)}U_{\e, y}, \sqrt{l(\bb)}U_{\e, y})$ is a positive solution of (\ref{eq4}). Note
that $2p=2^\ast$,
This implies that
\begin{align*}\lim_{\bb\to 0}\Big(k(\bb)+l(\bb)\Big)=k(0)+l(0)=\mu_1^{-\frac{N-2}{2}}+\mu_2^{-\frac{N-2}{2}},\end{align*}
that is, there exists $0<\bb_1\le\bb_2$, such that
$$k(\bb)+l(\bb)>\min\left\{\mu_1^{-\frac{N-2}{2}},\quad \mu_2^{-\frac{N-2}{2}}\right\},\quad \forall\,\,\bb\in (0, \bb_1).$$
Combining this with (\ref{ffc6}) and (\ref{eqqq47}), we have
$$I(U, V)=A'=A<I(\sqrt{k(\bb)}U_{\e, y}, \sqrt{l(\bb)}U_{\e, y}),\quad \forall \,\,\bb\in (0, \bb_1),$$
that is, $(\sqrt{k(\bb)}U_{\e, y}, \sqrt{l(\bb)}U_{\e, y})$ is different positive solution of (\ref{eq4}) with respect to $(U, V)$. This completes the
proof.\hfill$\square$\\

Before ending this section, we need to study the following properties of $(U, V)$ obtained in Theorem \ref{th2}.

\begin{proposition}\label{proposition} Assume that $\bb>0$. Let $(U, V)$ be a positive radially symmetric least energy solution of (\ref{eq4}) obtained in
Thoerem \ref{th2}. Then there exists $C>0$ such that
$$U(x)+V(x)\le C (1+|x|)^{2-N},\quad |\nabla U(x)|+|\nabla V(x)|\le C(1+|x|)^{1-N}.$$
\end{proposition}

\noindent{\bf Proof. } Define the Kelvin transformation:
$$U^*(x):=|x|^{2-N}U\left(\frac{x}{|x|^2}\right),\quad V^*(x):=|x|^{2-N}V\left(\frac{x}{|x|^2}\right).$$
Then $U^*, V^*\in D^{1, 2}(\RN)$ and $(U^*, V^*)$ satisfies the same system (\ref{eq4}). Then by a standard Brezis-Kato type argument (\cite{BK}),
we see that
$U^*, V^*\in L^{\iy}(\RN)$. Therefore, there exists $C>0$ such that
\be\label{ccc}U(x)+V(x)\le C |x|^{2-N}.\ee
On the other hand, note that $U, V$ are radially symmetric decreasing. We also have $U, V\in L^{\iy}(\RN)$, and so
$$U(x)+V(x)\le C (1+|x|)^{2-N}.$$
 Moreover, standard elliptic regularity theory implies that $U, V\in C^2(\RN)$. We write $U(|x|)=U(x)$ for convenience. Then
$$(r^{N-1}U_r)_r=-r^{N-1}(\mu_1 U^{2^\ast-1}+\bb U^{2^\ast/2-1}V^{2^\ast/2}),$$
and so for any $R\ge 1$, we see from (\ref{ccc}) that
{\allowdisplaybreaks
\begin{align*}
R^{N-1}|U_r(R)|&\le |U_r(1)|+\int_{1}^{R}r^{N-1}(\mu_1 U^{2^\ast-1}+\bb U^{2^\ast/2-1}V^{2^\ast/2})\,dr\\
&\le C+C\int_{1}^{+\iy}r^{N-1}r^{-N-2}\,dr\le C.
\end{align*}
}%
Therefore, it is easy to see that $|\nabla U(x)|\le C(1+|x|)^{1-N}$ for some $C>0$. Similarly, $|\nabla V(x)|\le C(1+|x|)^{1-N}$.\hfill$\square$

\vskip0.1in

\s{Proof of Theorem \ref{th1}}
\renewcommand{\theequation}{3.\arabic{equation}}

In this section, we assume that $-\la_1(\Omega)<\la_1\le\la_2<0$. Recalling the definition of $B$ in (\ref{eq10}),
since
$$\int_{\om}(|\nabla u|^2+\la_i u^2)\ge\left(1+\frac{\la_i}{\la_1(\om)}\right)\int_{\om}|\nabla u|^2, \quad i=1,2,$$
it is standard to see that $B>0$. As has been pointed out in Section 1, by \cite{BN}
the Brezis-Nirenberg problem (\ref{Brezis-Nirenberg})
$$-\Delta u +\la_i u=\mu_i |u|^{2^\ast-2}u,\qquad u\in H^1_0(\om)$$
has a positive least energy solution $u_{\mu_i}\in C^2(\om)\cap C(\overline{\om})$ with energy
{\allowdisplaybreaks
\begin{align}\label{eq11}\frac{1}{N}\left(\frac{\la_1(\om)+\la_i}{\la_1(\om)}\right)^{\frac{N}{2}}\mu_i^{-\frac{N-2}{2}} S^{N/2}\le B_{\mu_i}
&:=\frac{1}{2}\int_{\om}(|\nabla u_{\mu_i}|^2+\la_i u^2_{\mu_i})-
\frac{1}{2^\ast}\int_{\om}\mu_i u^{2^\ast}_{\mu_i}\nonumber\\
& <\frac{1}{N}\mu_i^{-\frac{N-2}{2}}S^{N/2},\quad i=1,2.\end{align}
}%
The next lemma is very important, where we need the assumption $\la_1,\la_2<0$.

\bl\label{lemma1} Let $\beta<0$, then
$$B<\min\left\{B_{\mu_1}+\frac{1}{N}\mu_2^{-\frac{N-2}{2}} S^{N/2},\,\,B_{\mu_2}+\frac{1}{N}\mu_1^{-\frac{N-2}{2}} S^{N/2},\,\,A\right\}.$$
\el

\noindent{\bf Proof. } The idea of this proof comes from \cite{CZ}, but some arguments are more delicate. Let $\beta<0$.
Let $t_0>0$ such that
\be\label{eq51}\frac{N}{2}B_{\mu_1}t^2-\frac{N}{4p} B_{\mu_1} t^{2p}+\frac{1}{N}(\mu_1/2)^{-\frac{N-2}{2}}S^{N/2}<0,\quad\forall \,\,t>t_0.\ee
Since $u_{\mu_1}\in C(\overline{\om})$ and $u_{\mu_1}\equiv 0$ on $\partial\Omega$, there exists
$B(y_0, 2R):=\{x : |x-y_0|\le 2R\}\subset \om$
such that
\be\label{eq52}\delta:=\max_{B(y_0, 2R)} u_{\mu_1}\le \min\left\{\left(\frac{\mu_2}{2|\beta|}\right)^{\frac{1}{p-1}},
\,\,\left(\frac{\la_1+\la_1(\om)}{2|\beta|}\right)^{
\frac{1}{p-1}}\right\}.\ee

Let $\psi\in C_0^1(B(y_0, 2R))$ be a function with $0\le\psi\le 1$ and $\psi\equiv 1$ for $|x-y_0|\le R$.
Define $v_\e=\psi U_{\e, y_0}$, where $U_{\e, y_0}$ is defined in (\ref{A-T}) and (\ref{A-T1}). Then by \cite{BN} or \cite[Lemma 1.46]{W},
we have the following inequalities
{\allowdisplaybreaks
\begin{gather}
   \label{eq53} \int_{\om}|\nabla v_\e|^2=S^{N/2}+O(\e^{N-2}),\quad\int_{\om}|v_\e|^{2^\ast}=S^{N/2}+O(\e^N),\\
    \label{eq54}\int_{\om}|v_\e|^2\ge C \e^2+O(\e^{N-2}).
\end{gather}
}%
Moreover, since $N\ge 5$, we have
{\allowdisplaybreaks
\begin{align}\label{eqq3}
\int_{\om}v_\e^{\frac{N}{N-2}}\,dx&\le\int_{B(y_0, 2R)}U_{\e, y_0}^{\frac{N}{N-2}}\,dx=C\int_{B(0, 2R)}\left(\frac{\e}{\e^2+|x|^2}\right)^{N/2}\,dx\nonumber\\
&\le C\e^{N/2}\left(\ln{\frac{2R}{\e}}+1\right)=o(\e^2).
\end{align}
}%
Since $supp(v_\e)\subset B(y_0, 2R)$, by (\ref{eq52}) we have for $t, s>0$ that
{\allowdisplaybreaks
\begin{align}\label{eq505}
2|\beta|t^ps^p\int_{\om}u_{\mu_1}^pv_\e^p &\le 2|\beta|\delta^{p-1} t^ps^p\int_{\om}u_{\mu_1}v_\e^p\nonumber\\
&\le|\beta|\delta^{p-1} t^{2p}\int_{\om}u_{\mu_1}^2
+|\beta|\delta^{p-1} s^{2p}\int_{\om}v_\e^{2p}\nonumber\\
&\le\frac{|\beta|\delta^{p-1}}{\la_1+\la_1(\om)}t^{2p}\int_{\om}(|\nabla u_{\mu_1}|^2+\la_1 u_{\mu_1}^2)+|\beta|\delta^{p-1}
s^{2p}\int_{\om}v_\e^{2p}\nonumber\\
&\le \frac{1}{2}t^{2p}\int_{\om}\mu_1 u_{\mu_1}^{2p}+\frac{1}{2} s^{2p}\int_{\om}\mu_2 v_\e^{2p},
\end{align}
}%
and so
{\allowdisplaybreaks
\begin{align}\label{eq55}
E(tu_{\mu_1}, sv_\e)&=\frac{1}{2}t^2\int_{\om}(|\nabla u_{\mu_1}|^2+\la_1 u_{\mu_1}^2)
   +\frac{1}{2}s^2\int_{\om}(|\nabla v_\e|^2+\la_2 v_\e^2)\nonumber\\
&\quad-\frac{1}{2p}\int_{\om}(t^{2p}\mu_1 u_{\mu_1}^{2p}+2t^ps^p\beta u_{\mu_1}^pv_\e^p +s^{2p}\mu_2 v_\e^{2p})\nonumber\\
&\le \frac{1}{2}t^2\int_{\om}(|\nabla u_{\mu_1}|^2 +\la_1 u_{\mu_1}^2)-\frac{1}{4p}t^{2p}\int_{\om}\mu_1u_{\mu_1}^{2p}\nonumber\\
&\quad+\frac{1}{2}s^2\int_{\om}(|\nabla v_{\e}|^2 +\la_2 v_{\e}^2)-\frac{1}{2p}s^{2p}\frac{\mu_2}{2}\int_{\om}v_{\e}^{2p}\nonumber\\
&=f(t)+g(s).
\end{align}
}%
By (\ref{eq53})-(\ref{eq54}), it is standard to check that (cf. \cite{BN, W})
\be\label{eq500}\max_{s>0}g(s)<\frac{1}{N}(\mu_2/2)^{-\frac{N-2}{2}} S^{N/2}\quad \hbox{for $\e$ small enough}.\ee
By (\ref{eq11}) we see that
$$f(t)=\frac{N}{2}B_{\mu_1}t^2-\frac{N}{4p} B_{\mu_1} t^{2p}.$$
Combining these with (\ref{eq51}), we get that
$$f(t)+ g(s)<0,\quad \forall\, t>t_0, \, \, s>0,$$
and so it follows from (\ref{eq55}) that
\begin{align*}
\max_{t, s>0}E(tu_{\mu_1}, sv_\e)=\max_{0<t\le t_0, s>0}E(tu_{\mu_1}, sv_\e).
\end{align*}
Define
$$
g_\e(s):=\frac{1}{2}s^2\int_{\om}(|\nabla v_\e|^2+\la_2 v_\e^2)\,dx-\frac{s^{2p}}{2p}\int_{\om}\mu_2v_\e^{2p}\,dx, \quad s>0.
$$
Then there exists a unique $s(\e)>0$, such that $g_\e'(s(\e))=0$ with
{\allowdisplaybreaks
\begin{align*}
s(\e)^{2p-2}&=\frac{\int_{\om}(|\nabla v_\e|^2+\la_2 v_\e^2)\,dx}{\int_{\om}\mu_2v_\e^{2p}\,dx}
\ge\left(1+\frac{\la_2}{\la_1(\om)}\right)\frac{\int_{\om}|\nabla v_\e|^2\,dx}{\mu_2\intR{U_{\e, y_0}^{2p}}}\\
&=\left(1+\frac{\la_2}{\la_1(\om)}\right)\frac{S^{N/2}+O(\e^{N-2})}{\mu_2 S^{N/2}}\\
&\ge\frac{1}{2\mu_2}\left(1+\frac{\la_2}{\la_1(\om)}\right)=: s_0^{2p-2}, \quad \hbox{for $\e$ small enough}.
\end{align*}
}%
Therefore, since $g_\e$ is increasing for $0<s\le s(\e)$, for any $0<s< s_0$, we have $g_\e(s)<g_\e(s_0)$ and so $E(t u_{\mu_1}, sv_\e)<E(t u_{\mu_1}, s_0v_\e)$.
That is,
\begin{align}\label{eq56}
\max_{t, s>0}E(tu_{\mu_1}, sv_\e)=\max_{0<t\le t_0, s\ge s_0}E(tu_{\mu_1}, sv_\e).
\end{align}
For $0<t\le t_0, s\ge s_0$, we see from (\ref{eqq3}) that
$$|\beta|t^ps^p\int_{\om}u_{\mu_1}^pv_\e^p\le |\beta|t_0^p\delta^p s_0^{p-2} s^2\int_{\om}v_\e^p\le C s^2\cdot o(\e^2),$$
and so
{\allowdisplaybreaks
\begin{align}\label{eq57}
E(tu_{\mu_1}, sv_\e)&=\frac{1}{2}t^2\int_{\om}(|\nabla u_{\mu_1}|^2+\la_1 u_{\mu_1}^2)
   +\frac{1}{2}s^2\int_{\om}(|\nabla v_\e|^2+\la_2 v_\e^2)\nonumber\\
&\quad-\frac{1}{2p}\int_{\om}(t^{2p}\mu_1 u_{\mu_1}^{2p}+2t^ps^p\beta u_{\mu_1}^pv_\e^p +s^{2p}\mu_2 v_\e^{2p})\nonumber\\
&\le \frac{1}{2}t^2\int_{\om}(|\nabla u_{\mu_1}|^2 +\la_1 u_{\mu_1}^2)-\frac{1}{2p}t^{2p}\int_{\om}\mu_1 u_{\mu_1}^{2p}\nonumber\\
&\quad+\frac{1}{2}s^2\left(\int_{\om}(|\nabla v_{\e}|^2 +\la_2 v_{\e}^2)+o(\e^2)\right)-\frac{1}{2p}s^{2p}\mu_2\int_{\om}v_{\e}^{2p}\nonumber\\
&=: f_1(t)+g_1(s).
\end{align}
}%
Note that $\max_{t>0}f_1(t)=f_1(1)=B_{\mu_1}$. By (\ref{eq53})-(\ref{eq54}) and $\la_2<0$, it is easy to show that
$$\max_{s>0}g_1(s)<\frac{1}{N}\mu_2^{-\frac{N-2}{2}} S^{N/2}\quad \hbox{for $\e$ small enough}.$$
Combining these with (\ref{eq56}) and (\ref{eq57}), we obtain that
{\allowdisplaybreaks
\begin{align}\label{eqq4}
\max_{t, s>0}E(tu_{\mu_1}, sv_\e)&=\max_{0<t\le t_0, s\ge s_0}E(tu_{\mu_1}, sv_\e)\nonumber\\
&\le\max_{t>0}f_1(t)+\max_{s>0}g_1(s)\nonumber\\
&<B_{\mu_1}+\frac{1}{N}\mu_2^{-\frac{N-2}{2}} S^{N/2}\quad \hbox{for $\e$ small enough}.
\end{align}
}%

Now, we claim that there exists $t_\e, s_\e>0$ such that $(t_\e u_{\mu_1}, s_\e v_\e)\in \mathcal{M}$.
Similarly as (\ref{eq505}), we have
{\allowdisplaybreaks
\begin{align*}
\left(\int_{\om}\bb u_{\mu_1}^pv_\e^p\,dx\right)^2&\le|\bb|^2\delta^{2p-2}\left(\int_{\om}u_{\mu_1}v_\e^p\,dx\right)^2\\
&\le|\bb|^2\delta^{2p-2}\int_{\om}u_{\mu_1}^2\,dx\int_{\om}v_\e^{2p}\,dx\\
&\le\frac{|\bb|^2\delta^{2p-2}}{(\la_1(\om)+\la_1)\mu_2}\int_{\om}\mu_1 u_{\mu_1}^{2p}\,dx\int_{\om}\mu_2 v_\e^{2p}\,dx\\
&<\int_{\om}\mu_1 u_{\mu_1}^{2p}\,dx\int_{\om}\mu_2 v_\e^{2p}\,dx.
\end{align*}
}%
For convenience we denote
\begin{gather*}
    D_1=\int_{\om}\mu_1 u_{\mu_1}^{2p}\,dx, \quad D_2=\int_{\om}\bb u_{\mu_1}^pv_\e^p\,dx,\\
    D_3=\int_{\om}\mu_2 v_\e^{2p}\,dx,\quad
D_4=\int_{\om}(|\nabla v_\e|^2+\la_2 v_\e^2)\,dx.
\end{gather*}
Then $D_2<0$ and $D_1D_3-D_2^2>0$. Furthermore, $(t u_{\mu_1}, sv_\e)\in\mathcal{M}$ for some $t, s>0$ is equivalent to
\be\label{eqq5}t^{2-p}D_1=t^p D_1+s^p D_2,\quad s^{2-p}D_4= s^p D_3+t^p D_2,\quad s, t>0.\ee
Note that $1<p=\frac{N}{N-2}<2$, by $s^p=(t^{2-p}-t^p)D_1/D_2>0$ we have $t>1$. Therefore, (\ref{eqq5}) is equivalent to
\be\label{eqq6}f_3(t):=D_4 \left(\frac{D_1}{|D_2|}(1-t^{2-2p})\right)^{\frac{2-p}{p}}-\frac{D_1D_3-D_2^2}{|D_2|}t^{2p-2}+\frac{D_1D_3}{|D_2|}, \,\,t>1.\ee
Since $f_3(1)>0$ and $\lim\limits_{t\to+\iy}f_3(t)<0$, (\ref{eqq6}) has a solution $t>1$. Hence (\ref{eqq5}) has
a solution $t_\e>0, s_\e>0$. That is, $(t_\e u_{\mu_1}, s_\e v_\e)\in \mathcal{M}$ and from (\ref{eqq4}) we get
$$B\le E(t_{\e}u_{\mu_1}, s_{\e}v_\e)\le\max_{t, s>0}E(t u_{\mu_1}, sv_\e)<B_{\mu_1}+\frac{1}{N}\mu_2^{-\frac{N-2}{2}} S^{N/2}.$$
By a similar argument, we can also prove that
$B<B_{\mu_2}+\frac{1}{N}\mu_1^{-\frac{N-2}{2}} S^{N/2}$. By (\ref{eq501}) and (\ref{eq11}), we have
$$A>\max\left\{B_{\mu_1}+\frac{1}{N}\mu_2^{-\frac{N-2}{2}} S^{N/2}, \,\,B_{\mu_2}+\frac{1}{N}\mu_1^{-\frac{N-2}{2}} S^{N/2}\right\},$$
This completes the proof.\hfill$\square$

\bl\label{lemma2}Assume that $\beta<0$, then there exists $C_2>C_1>0$,
such that for any $(u, v)\in \mathcal{M}$ with $E(u, v)\le A$, there holds
$$C_1\le\int_{\om} |u|^{2p}\,dx, \int_{\om}|v|^{2p}\,dx\le C_2.$$\el

\noindent{\bf Proof. } This follows directly from
{\allowdisplaybreaks
\begin{gather*}
    \frac{\la_1(\om)+\la_1}{\la_1(\om)}S\left(\int_{\om} |u|^{2p}\right)^{\frac{1}{p}}\le\int_{\om}(|\nabla u|^2+\la_1 u^2)
    \le\mu_1 \int_{\om}|u|^{2p},\\
    \frac{\la_1(\om)+\la_2}{\la_1(\om)}S\left(\int_{\om} |v|^{2p}\right)^{\frac{1}{p}}\le\int_{\om}(|\nabla v|^2+\la_2 v^2)
    \le\mu_1 \int_{\om}|v|^{2p},
\end{gather*}
}%
$E(u, v)\le A$ and (\ref{eq10}).\hfill$\square$

\bl\label{lem1}Let $u_n\rightharpoonup u, v_n\rightharpoonup v$ in $H^1_0(\om)$ as $n\to\iy$, then passing to a subsequence, there holds
$$\lim_{n\to\iy}\int_{\om}\left(|u_n|^p |v_n|^p-|u_n-u|^p |v_n-v|^p- |u|^p|v|^p\right)\,dx=0.$$\el

\noindent {\bf Proof. } Note that $2p=2^\ast$, we have
{\allowdisplaybreaks
\begin{gather*}
    u_n\to u, \,\,\,v_n\to v\quad\hbox{strongly in $L^{q}(\om), \,\,\,\forall\,\, 0< q<2p,$}\\
    u_n\rightharpoonup u,\,\,\, v_n\rightharpoonup v\quad\hbox{weakly in $L^{2p}(\om)$.}
\end{gather*}
}%
Fix any $t\in [0, 1]$. First, we claim that
\be\label{eqq7}|u_n-t u|^{p-2}(u_n-t u)|v_n|^p\rightharpoonup (1-t)^{p-1}|u|^{p-2}u|v|^p\,\,\hbox{weakly in $L^{\frac{2p}{2p-1}}(\om)$.}\ee
Since the map $h: L^{q_1}(\om)\to L^{q_1/q_2}(\om)$ with $h(s)=|s|^{q_2-1}s$ is continuous, so
{\allowdisplaybreaks
\begin{gather*}
    |u_n-t u|^{p-2}(u_n-tu)\to (1-t)^{p-1}|u|^{p-2}u\,\,\hbox{strongly in $L^{q}(\om), \,\,\,\forall\,\, 0< q<\frac{2p}{p-1},$}\\
    |v_n|^p\to |v|^p\,\,\hbox{strongly in $L^{q}(\om), \,\,\,\forall\,\, 0< q<2.$}
\end{gather*}
}%
Then for any $1\le q<\frac{2p}{2p-1}$, one has
$$|u_n-t u|^{p-2}(u_n-t u)|v_n|^p\to (1-t)^{p-1}|u|^{p-2}u|v|^p\quad\hbox{strongly in $L^{q}(\om)$.}$$
Since $|u_n-t u|^{p-2}(u_n-t u)|v_n|^p$ is uniformly bounded in $L^{\frac{2p}{2p-1}}(\om)$, passing to a subsequence, we may assume that
$|u_n-t u|^{p-2}(u_n-t u)|v_n|^p\rightharpoonup w$ weakly in $L^{\frac{2p}{2p-1}}(\om)$. Then for any $\vp\in C_0^{\iy}(\RN)$, we have
\begin{align*}
\int_{\om}w\vp=\lim_{n\to\iy}\int_{\om}|u_n-t u|^{p-2}(u_n-t u)|v_n|^p\vp=\int_{\om}(1-t)^{p-1}|u|^{p-2}u|v|^p\vp,
\end{align*}
which implies $w=(1-t)^{p-1}|u|^{p-2}u|v|^p$, that is, (\ref{eqq7}) holds. Similarly, we can show that
$|u_n-u|^{p}|v_n-t v|^{p-2}(v_n-t v)\rightharpoonup 0$ weakly in $L^{\frac{2p}{2p-1}}(\om)$.
Therefore, by (\ref{eqq7}), the Fubini Theorem and the Dominated Convergent Theorem,
{\allowdisplaybreaks
\begin{align*}
&\int_{\om}\left(|u_n|^p |v_n|^p-|u_n-u|^p |v_n-v|^p\right)\,dx\\
=&p\int_{\om}\int_0^1|u_n-t u|^{p-2}(u_n-tu) |v_n|^pu\,dt\,dx\\
&+p\int_{\om}\int_0^1|u_n- u|^{p} |v_n-t v|^{p-2}(v_n-tv)v\,dt\,dx\\
=&p\int_0^1\int_{\om}|u_n-t u|^{p-2}(u_n-tu) |v_n|^pu\,dx\,dt\\
&+p\int_0^1\int_{\om}|u_n- u|^{p} |v_n-t v|^{p-2}(v_n-tv)v\,dx\,dt\\
\to&p\int_0^1\int_{\om}(1-t)^{p-1}|u|^p|v|^p\,dx\,dt=\int_{\om}|u|^p|v|^p\,dx,\,\,\hbox{as $n\to\iy$.}
\end{align*}
}%
This completes the proof.\hfill$\square$\\

\noindent {\bf Proof of Theorem \ref{th1} for the case $\bb<0$. }
The main idea of the proof is similar to the proof of \cite[Theorem 1.3-(1)]{CZ} in case $N=4$,
but as we will see, some new ideas are needed. Assume that $\beta< 0$.
Note that $E$ is coercive and bounded from below on $\mathcal{M}$.
Then by the Ekeland variational priciple (cf. \cite{SM}), there exists a minimizing sequence $\{(u_n, v_n)\}\subset \mathcal{M}$ satisfying
{\allowdisplaybreaks
\begin{gather}
\label{eq20}E(u_n, v_n)\le \min\left\{B+\frac{1}{n}, \,\,A\right\},\\
\label{eq21}E(u, v)\ge E(u_n, v_n)-\frac{1}{n}\|(u_n, v_n)-(u, v)\|,\quad \forall (u, v)\in\mathcal{M}.
\end{gather}
}%
Here, $\|(u, v)\|:=(\int_{\om}(|\nabla u|^2+|\nabla v|^2)\,dx)^{1/2}$ is the norm of $H$. Then $\{(u_n, v_n)\}$ is bounded in $H$.
For any $(\va,\p)\in H$ with $\|\va\|,\|\p\|\le 1$ and each $n\in\mathbb{N}$, we define the functions $h_n$ and $g_n :\R^3\to \R$ by
{\allowdisplaybreaks
\begin{align}\label{eq22}
&h_n(t, s, l)= \int_{\om}|\nabla(u_n+t\va+s u_n)|^2+\la_1\int_{\om}|u_n+t\va+su_n|^2\nonumber\\
&\quad-\mu_1\int_{\om}|u_n+t\va+su_n|^{2p}-\beta\int_{\om}|u_n+t\va+su_n|^p|v_n+t\p+lv_n|^p,\\
\label{eq23}
&g_n(t, s, l)= \int_{\om}|\nabla(v_n+t\p+lv_n)|^2+\la_2\int_{\om}|v_n+t\p+lv_n|^2\nonumber\\
&\quad-\mu_2\int_{\om}|v_n+t\p+lv_n|^{2p}-\beta\int_{\om}|u_n+t\va+su_n|^p|v_n+t\p+lv_n|^p.
\end{align}
}%
Let $\mathbf{0}=(0,0,0)$. Then $h_n, g_n \in C^1(\R^3,\R)$, $h_n(\mathbf{0})=g_n(\mathbf{0})=0$ and
{\allowdisplaybreaks
\begin{gather*}
    \frac{\partial h_n}{\partial s}(\mathbf{0})=-(2p-2)\int_{\om}(|\nabla u_n|^2+\la_1 u_n^2)+p\bb\int_{\om}|u_n|^p|v_n|^p,\\
    \frac{\partial h_n}{\partial l}(\mathbf{0})=\frac{\partial g_n}{\partial s}(\mathbf{0})=-p\beta\int_{\om}|u_n|^p |v_n|^p\,dx,\\
    \frac{\partial g_n}{\partial l}(\mathbf{0})=-(2p-2)\int_{\om}(|\nabla v_n|^2+\la_2 v_n^2)+p\bb\int_{\om}|u_n|^p|v_n|^p.
\end{gather*}
}%
Define the matrix
$$
F_n:=\left(
     \begin{array}{cc}
       \frac{\partial h_n}{\partial s}(\mathbf{0}) & \frac{\partial h_n}{\partial l}(\mathbf{0}) \\
       \frac{\partial g_n}{\partial s}(\mathbf{0}) & \frac{\partial g_n}{\partial l}(\mathbf{0}) \\
     \end{array}
   \right).$$
Since $\beta<0$, it follows from Lemma \ref{lemma2} that
{\allowdisplaybreaks
\begin{align}\label{equ11}
\det(F_n)
&\ge (2p-2)^2\int_{\om}(|\nabla u_n|^2+\la_1|u_n|^2)\int_{\om}(|\nabla v_n|^2+\la_2|v_n|^2)\nonumber\\
&\ge CS^2\left(\int_{\om}|u_n|^{2p}\right)^{\frac{1}{p}}\left(\int_{\om}|v_n|^{2p}\right)^{\frac{1}{p}}\ge C>0,\end{align}
}%
where $C$ is independent of $n$.
By the implicit function theorem, functions $s_n(t)$ and $l_n(t)$
are well defined and class $C^1$ on some interval $(-\dd_n,+\dd_n)$ for $\dd_n>0$. Moreover, $s_n(0)=l_n(0)=0$
and
$$h_n(t, s_n(t), l_n(t))\equiv0,\quad g_n(t, s_n(t),l_n(t))\equiv0,\quad t\in(-\dd, +\dd).$$
With these, it is standard to prove that (see \cite[Theorem 1.3-(1)]{CZ} for instance)

\be\label{eq29}\lim_{n\to +\iy}E'(u_n, v_n)=0.\ee

Since $\{(u_n, v_n)\}$ is bounded in $H$, we may assume that $(u_n, v_n)\rightharpoonup (u, v)$ weakly in $H$. Passing to
a subsequence, we may assume that
{\allowdisplaybreaks
\begin{gather*}
    u_n\rightharpoonup u,\quad v_n\rightharpoonup v,\quad\hbox{weakly in}\,\,L^{2p}(\om),\\
    |u_n|^{q-1}u_n\rightharpoonup |u|^{q-1}u,\quad |v_n|^{q-1}v_n\rightharpoonup |v|^{q-1}v,\,\,\,\hbox{weakly in $L^{2p/q}(\om)$, $1< q<2p$,}\\
    u_n\to u,\quad v_n\to v,\quad\hbox{strongly in}\,\,L^2(\om).
\end{gather*}
}%
Thus, by (\ref{eq29}) we have $E'(u, v)=0$. Set $\omega_n=u_n-u$ and $\sigma_n=v_n-v$. Then by Brezis-Lieb Lemma (cf. \cite{W}),
there holds
\begin{gather}\label{eq30}
    |u_n|_{2p}^{2p}=|u|_{2p}^{2p}+|\omega_n|^{2p}_{2p}+o(1),\quad
    |v_n|_{2p}^{2p}=|v|_{2p}^{2p}+|\sigma_n|^{2p}_{2p}+o(1).
\end{gather}

Note that $(u_n, v_n)\in\mathcal{M}$ and $E'(u, v)=0$.
Combining these with (\ref{eq30}), Lemma \ref{lem1}, we get that
{\allowdisplaybreaks
\begin{gather}
\label{eq32}\int_{\om}|\nabla \omega_n|^2-\int_{\om}(\mu_1 |\omega_n|^{2p}+\beta |\omega_n|^p |\sigma_n|^p)=o(1),\\
\label{eq33}\int_{\om}|\nabla \sigma_n|^2-\int_{\om}(\mu_2 |\sigma_n|^{2p}+\beta |\omega_n|^p |\sigma_n|^p)=o(1),\\
\label{eq34}E(u_n, v_n)=E(u, v)+I(\omega_n, \sigma_n)+o(1).
\end{gather}
}%
Passing to a subsequence, we may assume that
$$\lim_{n\to+\iy}\int_{\om}|\nabla \omega_n|^2=b_1, \quad \lim_{n\to+\iy}\int_{\om}|\nabla \sigma_n|^2=b_2.$$
Then by (\ref{eq32}) and (\ref{eq33}) we have $I(\omega_n, \sigma_n)=\frac{1}{N}(b_1+b_2)+o(1)$. Letting $n\to +\iy$ in (\ref{eq34}),
we get that \be\label{eq35}0\le E(u, v)\le E(u, v)+ \frac{1}{N}(b_1+b_2)=\lim_{n\to+\iy}E(u_n, v_n)=B.\ee

{\bf Case 1.} $u\equiv 0, v\equiv 0$.

By Lemma \ref{lemma2}, (\ref{eq30}) and (\ref{eq35}), we have $0<b_1<+\iy$ and $0<b_2<+\iy$, and we may assume that
both $\omega_n\not\equiv 0$ and $\sigma_n\not\equiv 0$ for $n$ large.
Then by (\ref{eq32}) and (\ref{eq33}) we have
$$\int_{\om}\mu_1 |\omega_n|^{2p}\int_{\om}\mu_2 |\sigma_n|^{2p}-\left(\beta\int_{\om} |\omega_n|^p |\sigma_n|^p\right)^2>0,\quad\hbox{for $n$ large}.$$
Then by a similar argument as Lemma \ref{lemma1}, for $n$ large, there exists $t_n, s_n>0$ such that
$(t_n\omega_n, s_n\sigma_n)\in\mathcal{N}$. Up to a subsequence,
we claim that
\be\label{ffc1}\lim_{n\to+\iy}(|t_{n}-1|+|s_n-1|)=0.\ee
This conclusion is obvious in case $N=4$ and $p=2$ (see \cite{CZ}), but it is not trivial in our general case $N\ge 5$ here.  Denote
{\allowdisplaybreaks
\begin{gather*}
    B_{n, 1}=\int_{\om}|\nabla \omega_n|^2\to b_1,\,\,\, B_{n, 2}=\int_{\om}|\nabla \sigma_n|^2\to b_2,\\
    C_{n, 1}=\int_{\om}\mu_1 |\omega_n|^{2p},\,\,\,C_{n, 2}=\int_{\om}\mu_2 |\sigma_n|^{2p},\\
    D_n=|\beta|\int_{\om} |\omega_n|^p |\sigma_n|^p.
\end{gather*}
}%
Passing to a subsequence, we may assume that $C_{n, 1}\to c_1<+\iy$, $C_{n, 2}\to c_2<+\iy$ and $D_n\to d<+\iy$. By (\ref{eq32})-(\ref{eq33}) we have
\be\label{ffc0}c_1=b_1+d\ge b_1>0,\quad c_2=b_2+d\ge b_2>0,\ee
\be\label{ffff}t_n^2 B_{n, 1}= t_n^{2p} C_{n, 1}-t_n^p s_n^p D_n,\quad s_n^2 B_{n, 2}= s_n^{2p} C_{n, 2}-t_n^p s_n^p D_n.\ee
This implies that
\be\label{ffc}t_n^{2p-2}\ge \frac{B_{n, 1}}{C_{n, 1}}\to \frac{b_1}{c_1}>0,\quad s_n^{2p-2}\ge \frac{B_{n, 2}}{C_{n, 2}}\to \frac{b_2}{c_2}>0.\ee
Assume that, up to a subsequence, $t_n\to +\iy$ as $n\to \iy$, then by
$$t_n^{2p}C_{n,1}-t_n^2 B_{n, 1}=s_n^{2p}C_{n,2}-s_n^2 B_{n, 2},$$
we also have $s_n\to+\iy$.
Then
{\allowdisplaybreaks
\begin{align*}
d^2&=\lim_{n\to \iy}D_n^2=\lim_{n\to \iy}\frac{t_n^p C_{n, 1}-t_n^{2-p}B_{n, 1}}{s_n^p}\cdot\frac{s_n^p C_{n, 2}-s_n^{2-p}B_{n, 2}}{t_n^p}\\
&=\lim_{n\to \iy}(C_{n, 1}-t_n^{2-2p}B_{n, 1})(C_{n, 2}-s_n^{2-2p}B_{n, 2})\\
&=c_1c_2=(b_1+d)(b_2+d)>d^2,
\end{align*}
}%
a contradiction. Therefore, $t_n, s_n$ are uniformly bounded. Passing to a subsequence, by (\ref{ffc}) we may assume that $t_n\to t_\iy\ge
(b_1/c_1)^{\frac{1}{2p-2}}>0$ and $s_n\to s_\iy\ge (b_2/c_2)^{\frac{1}{2p-2}}>0$.
Then we see from (\ref{ffff}) that
$$s_\iy^p d= t_\iy^p c_1-t_\iy^{2-p}b_1,\quad t_\iy^p d = s_\iy^p c_2-s_\iy^{2-p}b_2.$$
If $d=0$, then $c_i=b_i$, and so $t_\iy=s_\iy=1$. That is, (\ref{ffc1}) holds. Now we consider the case $d>0$.
Define $f(t)=t^p c_1-t^{2-p}b_1$, then for $t\ge(b_1/c_1)^{\frac{1}{2p-2}}$, we have
\begin{align*}
f'(t)=pc_1 t^{p-1}-(2-p)b_1t^{1-p}>(2-p)t^{1-p}(c_1 t^{2p-2}-b_1)\ge 0,
\end{align*}
that is, $f$ is increasing with respect to $t\ge(b_1/c_1)^{\frac{1}{2p-2}}$. If $t_\iy<1$, then
$$s_\iy^p d=f(t_\iy)<f(1)=c_1-b_1=d,$$
that is, $s_\iy<1$, and we see from (\ref{ffc0}) that
{\allowdisplaybreaks
\begin{align*}
d^2&=\frac{t_\iy^p c_1-t_\iy^{2-p}b_1}{s_\iy^p}\cdot\frac{s_\iy^p c_2-s_\iy^{2-p}b_2}{t_\iy^p}=(c_1-t_\iy^{2-2p}b_1)(c_2-s_\iy^{2-2p}b_2)\\
&=(d+b_1-t_\iy^{2-2p}b_1)(d+b_2-s_\iy^{2-2p}b_2)<d^2,
\end{align*}
}%
a contradiction. If $t_\iy>1$, since $1\ge(b_1/c_1)^{\frac{1}{2p-2}}$, we have
$$s_\iy^p d=f(t_\iy)>f(1)=c_1-b_1=d,$$
that is, $s_\iy>1$, and so
\begin{align*}
d^2=(d+b_1-t_\iy^{2-2p}b_1)(d+b_2-s_\iy^{2-2p}b_2)>d^2,
\end{align*}
a contradiction.
Therefore, $t_\iy=s_\iy=1$ and (\ref{ffc1}) holds. This implies that
$$\frac{1}{N}(b_1+b_2)=\lim_{n\to+\iy}I(\omega_n, \sigma_n)=\lim_{n\to+\iy}I(t_n\omega_n, s_n\sigma_n)\ge A.$$
Combining this with (\ref{eq35}) we get that $B\ge A$, a contradiction with Lemma \ref{lemma1}. Therefore,
Case 1 is impossible.\\

{\bf Case 2.} $u\not\equiv 0, v\equiv 0$ or $u\equiv 0, v\not\equiv 0$.

Without loss of generality, we assume that $u\not\equiv 0, v\equiv 0$. Then $b_2>0$. By Case 1 we may assume that $b_1=0$. Then
$\lim_{n\to+\iy}\int_{\om} |\omega_n|^p |\sigma_n|^p=0$, and so
$$\int_{\om}|\nabla \sigma_n|^2=\int_{\om}\mu_2 |\sigma_n|^{2p}+o(1)\le \mu_2 S^{-p}\left(\int_{\om}|\nabla \sigma_n|^2\right)^p+o(1).$$
This implies that $b_2\ge \mu_2^{-\frac{N-2}{2}}S^{N/2}.$
Note that $u$ is a nontrivial solution of $-\Delta u+\la_1 u=\mu_1 |u|^{2^\ast-2}u$, we have from (\ref{eq11}) that $E(u, 0)\ge B_{\mu_1}$.
By (\ref{eq35}) we get that
$$B\ge B_{\mu_1}+\frac{1}{N}b_2\ge B_{\mu_1}+\frac{1}{N}\mu_2^{-\frac{N-2}{2}}S^{N/2},$$
a contradiction with Lemma \ref{lemma1}. Therefore,
Case 2 is impossible.\\

Since Cases 1 and 2 are both impossible, we have that $u\not\equiv 0, v\not\equiv 0$, that is, $(u, v)\in\mathcal{M}$. By (\ref{eq35}) we have
$E(u, v)=B$.
Then $(|u|, |v|)\in\mathcal{M}$ and $E(|u|, |v|)=B$. By Lemma \ref{lemma}, $(|u|, |v|)$ is a solution of (\ref{eq2}).
Then, using the maximum principle, we see that $|u|, |v|>0$ in $\om$. Therefore, $(|u|, |v|)$ is
a positive least energy solution of (\ref{eq2}). This completes the proof.\hfill$\square$\\

It remains to prove Theorem \ref{th1} for the case $\bb>0$. Let $\beta>0$ and define
\be\label{eq40} \mathcal{B}:=\inf_{h\in\Gamma}\max_{t\in[0,1]}E(h(t)),\ee
where $\Gamma = \{ h\in C([0,1], H) : h(0)=(0,0), E(h(1))<0\}.$ By (\ref{eq3}),
we see that for any $(u,v)\in H$, $(u,v)\neq(0,0)$,
{\allowdisplaybreaks
\begin{align}\label{eq41}
\max_{t>0}&E(t u, tv)=E(t_{u,v}u, t_{u,v}v)\nonumber\\
&=\frac{1}{N}t_{u,v}^2\int_{\om}(|\nabla u|^2+\la_1 u^2+|\nabla v|^2+\la_2 v^2)\nonumber\\
&=\frac{1}{N}t_{u,v}^{2^\ast}\int_{\om}(\mu_1 |u|^{2p}+2\beta |u|^p|v|^p+\mu_2 |v|^{2p}),
\end{align}
}%
 where $t_{u,v}>0$ satisfies \be\label{eq42}
t_{u,v}^{2p-2}=\frac{\int_{\om}(|\nabla u|^2+\la_1 u^2+|\nabla v|^2+\la_2 v^2)}{\int_{\om}(\mu_1 |u|^{2p}+2\beta |u|^p|v|^p+\mu_2 |v|^{2p})}.\ee
Note that $(t_{u,v}u, t_{u, v}v)\in\mathcal{M}'$, where
{\allowdisplaybreaks
\begin{align}\label{eq443}\mathcal{M}':=\Big\{(u, v)\in H\setminus\{ (0,0)\},&\,\,\, G(u, v):=\int_{\om}
(|\nabla u|^2+\la_1 u^2+|\nabla v|^2+\la_2 v^2)\nonumber\\
-& \int_{\om}(\mu_1 |u|^{2p}+2\beta |u|^p|v|^p+\mu_2 |v|^{2p})=0\Big\}, \end{align}
}%
it is easy to check that
\be\label{eq43} \mathcal{B}=\inf_{H\ni(u,v)\neq(0,0)}\max_{t>0}E(t
u,tv)=\inf_{(u, v)\in\mathcal{M}'}E(u, v).\ee
Note that $\mathcal{M}\subset\mathcal{M}'$, one has that
$\mathcal{B}\le B$. Similarly as (\ref{equation}), we have $\mathcal{B}>0$.

\bl\label{lemma3}Let $\bb>0$,
then
$$\mathcal{B}<\min\{B_{\mu_1}, \,\, B_{\mu_2},\,\, A\}.$$ \el

\noindent{\bf Proof. } Step 1. We prove that $\mathcal{B}< A$.
Without loss of generality, we may assume that $0\in \om$. Then there exists $\rho>0$
such that $B(0, 2\rho):=\{x: |x|\le 2\rho\}\subset\om$. Let $\psi\in C_0^1(B(0, 2\rho))$ be a nonnegative function with $0\le \psi\le 1$ and $\psi\equiv 1$ for
$|x|\le\rho$.
Recall that $(U, V)$ in Theorem \ref{th2}. we define
$$(U_\e(x), V_\e(x)):=\left(\e^{-\frac{N-2}{2}}U\left(\frac{x}{\e}\right), \,\,\e^{-\frac{N-2}{2}}V\left(\frac{x}{\e}\right)\right).$$
Then it is easy to see that
{\allowdisplaybreaks
\begin{align*}&\int_{\RN}|\nabla U_\e|^2=\int_{\RN}|\nabla U|^2,\quad \int_{\RN}|U_\e|^{2^\ast}=\int_{\RN}| U|^{2^\ast},\\
&\int_{\RN}|\nabla V_\e|^2=\int_{\RN}|\nabla V|^2,\quad \int_{\RN}|V_\e|^{2^\ast}=\int_{\RN}| V|^{2^\ast}.
\end{align*}
}%
Define
\be\label{equa6}(u_\e, v_\e):=(\psi U_\e, \psi V_\e).\ee
First we claim the following inequalities
{\allowdisplaybreaks
\begin{align}
    \label{cc1}&\int_{\om}|\nabla u_\e|^2\le\int_{\RN}|\nabla U|^2+O(\e^{N-2}),\\
   \label{cc2}&\int_{\om}|u_\e|^{2^\ast}\ge\int_{\RN}| U|^{2^\ast}+O(\e^N),\\
   \label{cc3}&\int_{\om}|u_\e|^{\frac{2^\ast}{2}}|v_\e|^{\frac{2^\ast}{2}}\ge\int_{\RN}| U|^{\frac{2^\ast}{2}}| V|^{\frac{2^\ast}{2}}+O(\e^N),\\
    \label{cc4}&\int_{\om}|u_\e|^2\ge C \e^2+O(\e^{N-2}),
\end{align}
}%
where $C$ is a positive constant.

Let $0<\e\ll\rho$. By Proposition \ref{proposition} we have
{\allowdisplaybreaks
\begin{align*}
\int_{\om}|\nabla \psi|^2|U_\e|^2\,dx&\le C\int_{\rho\le|x|\le2\rho}\e^{2-N}U^2(x/\e)\,dx\\
&\le C\e^2\int_{\rho/\e\le|x|\le2\rho/\e}U^2(x)\,dx\\
&\le C\e^2\int_{\rho/\e\le|x|\le2\rho/\e}|x|^{4-2N}\,dx=O(\e^{N-2});\\
\int_{\om}|\nabla U_\e|^2|\psi|^2\,dx&\le\int_{\RN}|\nabla U_\e|^2=\int_{\RN}|\nabla U|^2;\\
\left|\int_{\om}\psi U_\e\nabla \psi\nabla U_\e\right|&\le C\int_{\rho\le |x|\le 2\rho}|\nabla U_\e||U_\e|\,dx\\
&\le C\int_{\rho\le |x|\le 2\rho}\e^{1-N}|\nabla_x U(x/\e)||U(x/\e)|\,dx\\
&=C\e\int_{\rho/\e\le|x|\le2\rho/\e}|\nabla U(x)||U(x)|\,dx\\
&\le C\e \int_{\rho/\e\le|x|\le2\rho/\e}|x|^{3-2N}\,dx=O(\e^{N-2}).
\end{align*}
}%
Therefore,
{\allowdisplaybreaks
\begin{align*}
\int_{\om}|\nabla u_\e|^2\,dx&=\int_{\RN}|\nabla U_\e|^2|\psi|^2+\int_{\om}|\nabla \psi|^2|U_\e|^2+2\int_{\om}\psi U_\e\nabla \psi\nabla U_\e\\
&\le \int_{\RN}|\nabla U|^2\,dx+O(\e^{N-2}),
\end{align*}
}%
that is, (\ref{cc1}) holds. Note that
{\allowdisplaybreaks
\begin{align*}
\int_{\RN}(1-\psi^{2^\ast})|U_\e|^{2^\ast}\,dx&\le\int_{|x|\ge\rho}\e^{-N}|U(x/\e)|^{2^\ast}\,dx
=\int_{|x|\ge\rho/\e}|U(x)|^{2^\ast}\,dx\\
&\le C\int_{|x|\ge\rho/\e}|x|^{-2N}\,dx=O(\e^N),
\end{align*}
}%
then
{\allowdisplaybreaks
\begin{align*}
\int_{\om}|u_\e|^{2^\ast}\,dx&=\int_{\RN}|U_\e|^{2^\ast}\,dx-\int_{\RN}(1-\psi^{2^\ast})|U_\e|^{2^\ast}\,dx\\
&\ge\int_{\RN}|U|^{2^\ast}\,dx+O(\e^N),
\end{align*}
}%
that is, (\ref{cc2}) holds. Similarly, (\ref{cc3}) holds. Note that
{\allowdisplaybreaks
\begin{align*}
\int_{\om}|u_\e|^2\,dx&\ge\int_{|x|\le\rho}\e^{2-N}|U(x/\e)|^2\,dx\\
&=\e^2\int_{\RN}U^2\,dx-\e^2\int_{|x|\ge\rho/\e}U^2(x)\,dx\\
&\ge C\e^2-C\e^2\int_{|x|\ge\rho/\e}|x|^{4-2N}\,dx=C\e^2+O(\e^{N-2}),
\end{align*}
}%
that is, (\ref{cc4}) holds. Similarly, we have
{\allowdisplaybreaks
\begin{align}
    \label{cc9}&\int_{\om}|\nabla v_\e|^2\le\int_{\RN}|\nabla V|^2+O(\e^{N-2}),\\
   \label{cc10}&\int_{\om}|v_\e|^{2^\ast}\ge\int_{\RN}| V|^{2^\ast}+O(\e^N),\\
    \label{cc11}&\int_{\om}|v_\e|^2\ge C \e^2+O(\e^{N-2}).
\end{align}
}%

Recall that $I(U, V)=A$, we have
$$NA=\int_{\RN}|\nabla U|^2+|\nabla V|^2=\int_{\RN}\mu_1 U^{2^\ast}+2\bb U^{\frac{2^\ast}{2}}V^{\frac{2^\ast}{2}}+\mu_2V^{2^\ast}.$$
Combining this with (\ref{cc1})-(\ref{cc11}) and recalling that $\la_1,\la_2<0, 2p=2^\ast, N\ge 5$, we have for any $t>0$ that
{\allowdisplaybreaks
\begin{align}\label{eq13}
E(tu_\e, tv_\e)&=\frac{1}{2}t^2\int_{\om}(|\nabla u_\e|^2+\la_1 u_\e^2+|\nabla v_\e|^2+\la_2 v_\e^2)
   \nonumber\\
&\quad-\frac{1}{2p}t^{2p}\int_{\om}(\mu_1 u_\e^{2p}+2\beta u_\e^pv_\e^p +\mu_2 v_\e^{2p})\nonumber\\
&\le \frac{1}{2}\left(\int_{\RN}\left(|\nabla U|^2+|\nabla V|^2\right)- C\e^2+ O(\e^{N-2})\right)t^2\nonumber\\
&\quad-\frac{1}{2^\ast}\left(\int_{\RN}\left(\mu_1 U^{2^\ast}+2\bb U^{\frac{2^\ast}{2}}V^{\frac{2^\ast}{2}}+\mu_2V^{2^\ast}\right)
+O(\e^N)\right)t^{2^\ast}\nonumber\\
&= \frac{1}{2}\left(NA- C\e^2+ O(\e^{N-2})\right)t^2
-\frac{1}{2^\ast}\left(NA +O(\e^N)\right)t^{2^\ast}\nonumber\\
&\le \frac{1}{N}\Big(NA- C\e^2+ O(\e^{N-2})\Big)\left(\frac{NA- C\e^2+ O(\e^{N-2})}{NA +O(\e^N)}\right)^{\frac{N-2}{2}}\nonumber\\
&<A\quad\hbox{for $\e>0$ small enough.}
\end{align}
}%
Hence, for $\e>0$ small enough, there holds
\begin{align}\label{equa7}
\mathcal{B}\le\max_{t>0}E(tu_\e, tv_\e)<A.
\end{align}

Step 2. we shall prove that $\mathcal{B}<B_{\mu_1}$. This proof is similar to Lemma \ref{lemm2}. Recall (\ref{eq11}) and (\ref{eq42}), we define $t(s):=
t_{u_{\mu_1}, su_{\mu_2}}$, that is,
$$t(s)^{2p-2}=\frac{NB_{\mu_1}+s^2 NB_{\mu_2}}{NB_{\mu_1}+|s|^{2p} NB_{\mu_2}+|s|^p\int_{\om}2\beta |u_{\mu_1}|^p|u_{\mu_2}|^p}.$$
Note that $t(0)=1$. Recall that $1<p=\frac{N}{N-2}<2$, by direct computations we get that
$$\lim_{s\to 0}\frac{t'(s)}{|s|^{p-2}s}=-\frac{p\int_{\om}2\beta |u_{\mu_1}|^p|u_{\mu_2}|^p}{(2p-2)NB_{\mu_1}},$$
that is,
$$t'(s)=-\frac{p\int_{\om}2\beta |u_{\mu_1}|^p|u_{\mu_2}|^p}{(2p-2)NB_{\mu_1}}|s|^{p-2}s(1+o(1)),\,\,\,\hbox{as $s\to 0$,}$$
and so
$$t(s)=1-\frac{\int_{\om}2\beta |u_{\mu_1}|^p|u_{\mu_2}|^p}{(2p-2)NB_{\mu_1}}|s|^{p}(1+o(1)),\,\,\,\hbox{as $s\to 0$.}$$
This implies that
$$t(s)^{2p}=1-\frac{2p\int_{\om}2\beta |u_{\mu_1}|^p|u_{\mu_2}|^p}{(2p-2)NB_{\mu_1}}|s|^{p}(1+o(1)),\,\,\,\hbox{as $s\to 0$.}$$
Therefore, we deduce from (\ref{eq41}) and $\frac{2p}{2p-2}=N/2$ that
{\allowdisplaybreaks
\begin{align*}
\mathcal{B}&\le E\left(t(s)u_{\mu_1}, t(s)su_{\mu_2}\right)\\
&=\frac{t(s)^{2p}}{N}\left(NB_{\mu_1}+|s|^{2p} NB_{\mu_2}+|s|^p\int_{\om}2\beta |u_{\mu_1}|^p|u_{\mu_2}|^p\right)\\
&=B_{\mu_1}-\left(\frac{1}{2}-\frac{1}{N}\right)|s|^p\int_{\om}2\beta |u_{\mu_1}|^p|u_{\mu_2}|^p+o(|s|^p)\\
&<B_{\mu_1}\quad \hbox{as $|s|>0$ small enough,}
\end{align*}
}%
that is, $\mathcal{B}<B_{\mu_1}$. By a similar argument, we can prove that $\mathcal{B}<B_{\mu_2}$. This completes the proof.
\hfill$\square$\\

\noindent {\bf Proof of Theorem \ref{th1} for the case $\bb>0$. } Assume that $\bb>0$. Since the functional $E$ has a mountain pass structure,
by the mountain pass theorem (cf. \cite{AR-pass, W}) there exists $\{(u_n, v_n)\}\subset H$ such that
$$\lim_{n\to+\iy}E(u_n, v_n)=\mathcal{B},\quad \lim_{n\to+\iy}E'(u_n, v_n)=0.$$
It is standard to see that $\{(u_n, v_n)\}$ is bounded in $H$, and so we may assume that $(u_n, v_n)\rightharpoonup (u, v)$ weakly in $H$.
Set $\omega_n=u_n-u$ and $\sigma_n=v_n-v$ and use the same symbols as in the proof of Theorem \ref{th1} for the case $\bb<0$,
we see that $E'(u, v)=0$ and (\ref{eq32})-(\ref{eq34}) also hold. Moreover,
\be\label{eq46}0\le E(u, v)\le E(u, v)+ \frac{1}{N}(b_1+b_2)=\lim_{n\to+\iy}E(u_n, v_n)=\mathcal{B}.\ee

{\bf Case 1.} $u\equiv 0, v\equiv 0$.

By (\ref{eq46}), we have $b_1+b_2>0$. Then we may assume that $(\omega_n, \sigma_n)\neq (0, 0)$ for $n$ large.
Recall $\mathcal{N}'$ in (\ref{eqqq41}), by (\ref{eq32})-(\ref{eq33}), it is
easy to check that there exists $t_n>0$ such that $(t_n \omega_n, t_n \sigma_n)\in \mathcal{N}'$ and $t_n\to 1$ as $n\to \iy$.
Then by (\ref{ffc2}) and (\ref{eq46}) we have
$$\mathcal{B}=\frac{1}{N}(b_1+b_2)=\lim_{n\to+\iy}I(\omega_n, \sigma_n)=\lim_{n\to+\iy}I(t_n\omega_n, t_n\sigma_n)\ge A'=A,$$
a contradiction with Lemma \ref{lemma3}. Therefore,
Case 1 is impossible.\\

{\bf Case 2.} $u\not\equiv 0, v\equiv 0$ or $u\equiv 0, v\not\equiv 0$.

Without loss of generality, we may assume that $u\not\equiv 0, v\equiv 0$. Then $u$ is a nontrivial
solution of $-\Delta u+\la_1 u=\mu_1 |u|^{2^\ast-2}u$, and so $\mathcal{B}\ge E(u, 0)\ge B_{\mu_1}$,
a contradiction with Lemma \ref{lemma3}. Therefore,
Case 2 is also impossible.\\

Since Cases 1 and 2 are both impossible, we have that $u\not\equiv 0, v\not\equiv 0$. Since $E'(u, v)=0$, we have
$(u, v)\in\mathcal{M}$. By $\mathcal{B}\le B$ and (\ref{eq46}) we have
$E(u, v)=\mathcal{B}=B$. This means $(|u|, |v|)\in \mathcal{M}\subset\mathcal{M}'$ and
$E(|u|, |v|)=\mathcal{B}=B$. By (\ref{eq443}) and (\ref{eq43}), there exists a Lagrange multiplier
$\gamma\in \R$ such that
$$E'(|u|, |v|)-\gamma G'(|u|, |v|)=0.$$
Since $E'(|u|, |v|)(|u|, |v|)=G(|u|, |v|)=0$ and $$G'(|u|, |v|)(|u|, |v|)=-(2p-2)\int_{\om}(\mu_1 |u|^{2p}+2\beta |u|^p|v|^p+\mu_2 |v|^{2p})\neq0,$$
we get that $\gamma=0$ and so $E'(|u|, |v|)=0$.
This means that
$(|u|, |v|)$ is a least energy solution of (\ref{eq2}). By the maximum principle, we see that $|u|, |v|>0$ in $\om$. Therefore, $(|u|, |v|)$ is
a positive least energy solution of (\ref{eq2}). \hfill$\square$

\vskip0.1in

\s{Proof of Theorems \ref{unique1} and \ref{unique2}}
\renewcommand{\theequation}{4.\arabic{equation}}

In this section, we assume that $-\la_1(\om)<\la_1= \la_2=\la<0$ and $\bb\ge (p-1)\max\{\mu_1, \mu_2\}$.
Define $g:[(p-1)\max\{\mu_1, \mu_2\},\,\,+\iy)$ by
\be\label{eqq001}g(\bb):=(p-1)\mu_1\mu_2\bb^{2/p-2}+\bb^{2/p}.\ee
Then
$$g'(\bb)=\frac{2}{p}\bb^{2/p-3}\Big(\bb^2-(p-1)^2\mu_1\mu_2\Big)>0,\quad \forall\,\,\bb>(p-1)\max\{\mu_1, \mu_2\}.$$
By direct computations, we have
$$
g\Big((p-1)\max\{\mu_1, \mu_2\}\Big)\le p(p-1)^{\frac{2}{p}-1}\max\left\{\mu_1^{2/p},\,\,\mu_2^{2/p}\right\}.
$$
Therefore, there exists a unique $\bb_0\ge (p-1)\max\{\mu_1, \mu_2\}$ such that
\be\label{eqq002}g(\bb_0)=p(p-1)^{\frac{2}{p}-1}\max\left\{\mu_1^{2/p},\,\,\mu_2^{2/p}\right\},\quad\hbox{and}\ee
\be\label{eqq003}g(\bb)>p(p-1)^{\frac{2}{p}-1}\max\left\{\mu_1^{2/p},\,\,\mu_2^{2/p}\right\},\quad\forall\,\,\bb>\bb_0.\ee
Moreover,
\be\label{eqq004}\bb_0=(p-1)\max\{\mu_1, \mu_2\},\quad\hbox{if $\mu_1=\mu_2$}.\ee

\bl\label{lem6} Assume that $\bb>\bb_0$, where $\bb_0$ is defined in (\ref{eqq002}). Let $(k_0, l_0)$ be in Lemma \ref{lem}. Then $p\mu_1 k_0^{p-1}<1$ and $p\mu_2 l_0^{p-1}<1$.\el

\noindent {\bf Proof. } Let $k_1=(p\mu_1)^{\frac{1}{p-1}}$, then by (\ref{eqq03}) we have
$$l_1:=h_1(k_1)=\left[\frac{p-1}{p\bb(p\mu_1)^{\frac{2-p}{2(p-1)}}}\right]^{2/p}.$$
By (\ref{eqq003}) and direct computations, we get that
{\allowdisplaybreaks
\begin{align*}
\al_2(k_1, l_1)&=\mu_2l_1^{p-1}+\bb k_1^{p/2}l_1^{p/2-1}-1\\
&=\frac{1}{l_1}\left[\mu_2 l_1^p+k_1(1-\mu_1 k_1^{p-1})\right]-1
=\frac{1}{l_1}\left[\mu_2 l_1^p+\frac{p-1}{p}k_1\right]-1\\
&=\left[\frac{p\bb(p\mu_1)^{\frac{2-p}{2(p-1)}}}{p-1}\right]^{\frac{2}{p}}\left\{\mu_2
\left[\frac{p-1}{p\bb(p\mu_1)^{\frac{2-p}{2(p-1)}}}\right]^{2}+\frac{p-1}{p}(p\mu_1)^{-\frac{1}{p-1}}\right\}-1\\
&=(p-1)^{1-2/p}p^{-1}\mu_1^{-2/p}g(\bb)-1>0.
\end{align*}
}%
Combining this with Lemma \ref{lem4} we have $k_1>k_0$, that is, $p\mu_1 k_0^{p-1}<1$. Similarly, let $l_2=(p\mu_2)^{\frac{1}{p-1}}$, then
$$\al_1(h_2(l_2), l_2)=(p-1)^{1-2/p}p^{-1}\mu_2^{-2/p}g(\bb)-1>0.$$
By Lemma \ref{lem4} again, we have $l_2>l_0$, and so $p\mu_2 l_0^{p-1}<1$.\hfill$\square$\\

\bl\label{lem7}Assume that $\bb>\bb_0$, where $\bb_0$ is defined in (\ref{eqq002}). Recall $\al_1, \al_2$ defined
in (\ref{eqq01})-(\ref{eqq02}), and $(k_0, l_0)$ obtained in Lemma \ref{lem}. Then
$$F(k_0, l_0):=\det\left(
                     \begin{array}{cc}
                       \partial_k \al_1(k_0, l_0) & \partial_l \al_1(k_0, l_0) \\
                       \partial_k \al_2(k_0, l_0) & \partial_l \al_2(k_0, l_0) \\
                     \end{array}
                   \right)<0.
$$\el

\noindent {\bf Proof. } By $\al_1(k_0, l_0)=\al_2(k_0, l_0)=0$ we have
$$\bb k_0^{p/2-2}l_0^{p/2}=k_0^{-1}-\mu_1 k_0^{p-2},\quad \bb l_0^{p/2-2}k_0^{p/2}=l_0^{-1}-\mu_2 l_0^{p-2}.$$
Then
{\allowdisplaybreaks
\begin{align*}
    \partial_k \al_1(k_0, l_0)&=(p-1)\mu_1 k_0^{p-2}+(p/2-1)\bb k_0^{p/2-2}l_0^{p/2}\\
    &=\frac{p}{2}\mu_1 k_0^{p-2}-(1-p/2)k_0^{-1};\\
   \partial_l \al_2(k_0, l_0)&=(p-1)\mu_2 l_0^{p-2}+(p/2-1)\bb l_0^{p/2-2}k_0^{p/2}\\
   &=\frac{p}{2}\mu_2 l_0^{p-2}-(1-p/2)l_0^{-1};\\
   \partial_l \al_1(k_0, l_0)&=\partial_k \al_2(k_0, l_0)=\frac{p}{2}\bb k_0^{p/2-1}l_0^{p/2-1}\\
   &=\frac{p}{2}\sqrt{(k_0^{-1}-\mu_1 k_0^{p-2})(l_0^{-1}-\mu_2 l_0^{p-2})}.
\end{align*}
}%
Therefore,
{\allowdisplaybreaks
\begin{align*}
F(k_0, l_0)=&\left[\frac{p}{2}\mu_1 k_0^{p-2}-(1-p/2)k_0^{-1}\right]\left[\frac{p}{2}\mu_2 l_0^{p-2}-(1-p/2)l_0^{-1}\right]\\
&-\frac{p^2}{4}\Big(k_0^{-1}-\mu_1 k_0^{p-2}\Big)\Big(l_0^{-1}-\mu_2 l_0^{p-2}\Big)\\
=&\frac{p}{2}(p-1)k_0^{-1}l_0^{-1}\left(\mu_1 k_0^{p-1}+\mu_2 l_0^{p-1}-\frac{2}{p}\right)<0
\end{align*}
}%
from Lemma \ref{lem6}.
\hfill$\square$

\bl\label{lem8} Fix any $\mu_1, \mu_2>0$ and $\bb>\bb_0$. Let $(u_0, v_0)$ be a least energy solution of (\ref{eq2})
with $(\mu_1, \mu_2, \bb)$ which exists by Theorem \ref{th1}. Recall $(\sqrt{k_0}\omega, \sqrt{l_0} \omega)$ in Theorem \ref{th0}. Then
\be\label{eqq008}\int_{\om}|u_0|^{2p}\,dx=k_0^p\int_{\om}\omega^{2p}\,dx.\ee\el

\noindent {\bf Proof. } Fix any $\mu_1, \mu_2>0$ and $\bb>\bb_0$. We remark from (\ref{eqq001})-(\ref{eqq002}) that $\bb_0(\mu_1, \mu_2):=\bb_0$ is
completely determined by $\mu_1,\mu_2$. Hence there exists $0<\e<\mu_1$ such
that for any $\mu\in (\mu_1-\e, \mu_1+\e)$, one has $\bb>\bb_0(\mu, \mu_2)$. Then by Lemmas \ref{lem}, \ref{lem7} and the implicit function theorem,
when $\mu_1$ is replaced by $\mu$, functions $k_0(\mu)$ and $l_0(\mu)$ are well defined and class $C^1$ for $\mu\in(\mu_1-\e_1, \mu_1+\e_1)$ for some $0<\e_1\le
\e$. Recall the definition of
$E, \mathcal{M}$ and $B$, they all depend on $\mu$, and we use notations $E_\mu, \mathcal{M}_\mu, B(\mu)$ in this proof, when $\mu_1$ is replaced by $\mu$.
Then $B(\mu)=(k_0(\mu)+l_0(\mu))B_1\in C^1((\mu_1-\e_1, \mu_1+\e_1), \R)$. In particular, $B'(\mu_1):=\frac{d}{d\mu}B(\mu_1)$ exists.
Note that $B=\mathcal{B}$ by the proof of Theorem \ref{th1} for the case $\bb>0$. Then by (\ref{eq43}) we have
$$B(\mu)=\inf_{H\ni(u,v)\neq(0,0)}\max_{t>0}E_\mu(t
u,tv).$$
Denote
{\allowdisplaybreaks
\begin{gather*}
  C=\int_{\om}(|\nabla u_0|^2+\la_1 u_0^2+|\nabla v_0|^2+\la_2 v_0^2),\\
  D=\int_{\om}(2\beta |u_0|^p|v_0|^p+\mu_2 |v_0|^{2p}),\quad G=\int_{\om}|u_0|^{2p}\,dx.
\end{gather*}
}%
There exists $t(\mu)>0$ such that
$$
\max_{t>0}E_\mu(t
u_0,tv_0)=E_{\mu}\Big(t(\mu) u_0, t(\mu)v_0\Big),
$$
where $t(\mu)>0$
satisfies
$f(\mu, t(\mu))=0$, and
$$f(\mu, t):=t^{2p-2}(\mu G+D)-C.$$
Note that $f(\mu_1, 1)=0$, $\frac{\partial}{\partial t}f(\mu_1, 1)=(2p-2)(\mu_1 G+ D)>0$, and $f(\mu, t(\mu))\equiv 0$. By the implicit function
theorem, there exists $0<\e_2\le \e_1$, such that $t(\mu)\in C^{\iy}((\mu_1-\e_2, \mu_1+\e_2), \R)$. By $f(\mu, t(\mu))\equiv 0$ we see that
$$t'(\mu_1)=-\frac{G}{(2p-2)(\mu_1 G+D)}.$$
By Taylor expansion one has
$t(\mu)=1+t'(\mu_1)(\mu-\mu_1)+O((\mu-\mu_1)^2)$, and so
\begin{align*}t^2(\mu)&=1+2t'(\mu_1)(\mu-\mu_1)+O((\mu-\mu_1)^2).
\end{align*}
Note that $C=\mu_1 G +D=N B(\mu_1)$. Then by (\ref{eq41}) that
{\allowdisplaybreaks
\begin{align*}
B(\mu)&\le E_{\mu}(t(\mu)u_0, t(\mu)v_0)=\frac{1}{N}t^2(\mu)C=t^2(\mu)B(\mu_1)\\
&=B(\mu_1)-\frac{2G B(\mu_1)}{(2p-2)(\mu_1 G+D)}(\mu-\mu_1)+O((\mu-\mu_1)^2)\\
&=B(\mu_1)-\frac{G}{2p}(\mu-\mu_1)+O((\mu-\mu_1)^2),
\end{align*}
}%
It follows that
$ \displaystyle  \frac{B(\mu) -B(\mu_1)}{\mu -\mu_1}\ge -\frac{G}{2p} +O ((\mu -\mu_1)),\,\,\hbox{ as}\,\, \mu\nearrow \mu_1,$
and so $B'(\mu_1) \ge -\frac{G}{2p}$. Similarly, we have
$$\frac{B(\mu) -B(\mu_1)}{\mu -\mu_1}\le -\frac{G}{2p} +O ((\mu -\mu_1)),\,\, \hbox{ as}\,\, \mu\searrow \mu_1,$$
that is, $B'(\mu_1) \le -\frac{G}{2p}$. Hence, $ \displaystyle B'(\mu_1) =
-\frac{G}{2p}=-\frac{1}{2p}\int_{\om}|u_0|^{2p}\,dx.$ By Theorem \ref{th0}, $(\sqrt{k_0}\omega, \sqrt{l_0}\omega)$ is also a positive least energy solution
of (\ref{eq2}). Therefore, $ \displaystyle B'(\mu_1)
=-\frac{k_0^p}{2p}\int_{\om}\omega^{2p}\,dx$, that is, (\ref{eqq008}) holds.\hfill$\square$\\

\noindent {\bf Proof of Theorem \ref{unique1}. } Let $(u, v)$ be any a positive least energy solution
of (\ref{eq2}). By Lemma \ref{lem8}, we have
$$\int_{\om}|u|^{2p}\,dx=k_0^p\int_{\om}\omega^{2p}\,dx.$$
By a similar proof of Lemma \ref{lem8}, that is, by computing $B'(\mu_2)$ and $B'(\bb)$ respectively, we can show that
\begin{gather*}
    \int_{\om}|v|^{2p}\,dx=l_0^p\int_{\om}\omega^{2p}\,dx,
    \quad\hbox{and}\,\,\,\int_{\om}|u|^p|v|^p\,dx=k_0^{p/2}l_0^{p/2}\int_{\om}\omega^{2p}\,dx.
\end{gather*}
Therefore,
\be\label{eqq009}\int_{\om}|u|^p|v|^p\,dx=l_0^{p/2}k_0^{-p/2}\int_{\om}|u|^{2p}\,dx,\quad \int_{\om}|u|^p|v|^p\,dx=l_0^{-p/2}k_0^{p/2}\int_{\om}|v|^{2p}\,dx.\ee
Define $(\tilde{u}, \tilde{v}):=(\frac{1}{\sqrt{k_0}}u, \frac{1}{\sqrt{l_0}}v).$ By $\al_1(k_0, l_0)=\al_2(k_0, l_0)=0$ and (\ref{eqq009}) we get
\be\label{eqq0010}\int_{\om}|\nabla \tilde{u}|^2+\la \tilde{u}^2\,dx=\int_{\om}|\tilde{u}|^{2p}\,dx,\quad \int_{\om}|\nabla \tilde{v}|^2+\la
\tilde{v}^2\,dx=\int_{\om}|\tilde{v}|^{2p}\,dx.\ee
Then by (\ref{func1}) we have
$$\frac{1}{N}\int_{\om}|\nabla \tilde{u}|^2+\la \tilde{u}^2\,dx\ge B_1,\quad \frac{1}{N}\int_{\om}|\nabla \tilde{v}|^2+\la \tilde{v}^2\,dx\ge B_1,$$
and so
{\allowdisplaybreaks
\begin{align*}
B=(k_0+l_0)B_1&=\frac{1}{N}\int_{\om}(|\nabla u|^2+\la_1 u^2+|\nabla v|^2+\la_2 v^2)\\
&=\frac{1}{N}k_0\int_{\om}(|\nabla \tilde{u}|^2+\la_1 \tilde{u}^2)+\frac{1}{N}l_0\int_{\om}(|\nabla \tilde{v}|^2+\la_1 \tilde{v}^2)\\
&\ge (k_0+l_0)B_1.
\end{align*}
}%
This implies that
$$\frac{1}{N}\int_{\om}|\nabla \tilde{u}|^2+\la \tilde{u}^2\,dx= B_1,\quad \frac{1}{N}\int_{\om}|\nabla \tilde{v}|^2+\la \tilde{v}^2\,dx= B_1.$$
Combining this with (\ref{eqq0010}), we see from \cite{BN} that $\tilde{u}$
and $\tilde{v}$ are both postive least energy solutions of (\ref{BN}).
Then we see from $(u, v)$ satisfies (\ref{eq2}) that
\begin{align*}
-\Delta \tilde{u}+\la \tilde{u}&=\mu_1 k_0^{p-1}\tilde{u}^{2p-1}+\bb k_0^{p/2-1}l_0^{p/2}\tilde{u}^{p-1}\tilde{v}^p=\tilde{u}^{2p-1},
\end{align*}
that is, $\tilde{u}^{p-1}\tilde{v}^p= \tilde{u}^{2p-1}$ and so $\tilde{u}=\tilde{v}$. Denote $U=\tilde{u}$, then $(u, v)=(\sqrt{k_0}U, \sqrt{l_0} U)$,
where $U$ is a positive least energy solution of (\ref{BN}).

Now we assume that $\om$ is a ball in $\RN$, then the positive least energy
solution of the Brezis-Nirenberg problem (\ref{BN}) is unique (cf. \cite{AY}).
Therefore,
the positive least energy solution of (\ref{eq2}) is unique.
\hfill$\square$\\

\noindent {\bf Proof of Theorem \ref{unique2}. } The proof is the same as that of Theorem \ref{unique1}.\hfill$\square$

\vskip0.1in

\s{Proof of Theorems \ref{th3} and \ref{sign-changing}}
\renewcommand{\theequation}{5.\arabic{equation}}

This section is devoted to the proof of Theorems \ref{th3} and \ref{sign-changing}. Recall the definition of $E, \mathcal {M}$ and $B$,
they both depend on $\beta$, and we use
notations $E_\beta, \mathcal{M}_\bb, B_\bb$ in this section. Define $B(x_0, R):=\{x\in\RN : |x-x_0|<R \}$.
Consider the problem
\be\label{eqq11}
\begin{cases}-\Delta u +\lambda_2 u= \mu_2 u^{2^{\ast}-1}\,\,\,\hbox{in} \,\,B(0,R),\\
u>0\,\,\,\hbox{in}\,\, B(0,R),\quad u=0 \,\,\,\hbox{on} \,\, \partial
B(0, R),\end{cases}\ee
and the corresponding functional is $J_{R}: H_0^1
(B(0, R)) \to \mathbb{R}$ given by
\be\label{eqq12} J_{R} (u)
= \frac{1}{2}\int_{B(0, R)}(|\nabla u|^2 + \lambda_2 u^2) \,dx
-\frac{1}{2^{\ast}}\mu_2\int_{B(0, R)}|u|^{2^{\ast}}\,dx.\ee
We need the following energy estimates from the authors' paper \cite{CZ1}.

\bt\label{th-5} (see \cite{CZ1}) Let $N\ge 5$. Then there exists $R_0>0$ and $C_1, C_2>0$, such that for any $0<R<R_0$, (\ref{eqq11}) has a
least energy solution $U_R$ and
\be\label{eqq13}\frac{1}{N}\mu_2^{-\frac{N-2}{2}}S^{N/2}- C_1 R^{\frac{2N-4}{N-4}}\le
J_{R}(U_R)\le \frac{1}{N}\mu_2^{-\frac{N-2}{2}}S^{N/2}- C_2
R^{\frac{2N-4}{N-4}}.\ee\et

With the help of Theorem \ref{th-5}, we have the following lemma, which improves Lemma \ref{lemma1} in case $N\ge 6$.

\bl\label{lem2}Let $N\ge 6$. Then
$$\sup\limits_{\bb<0}B_\bb<\min\left\{B_{\mu_1}+\frac{1}{N}\mu_2^{-\frac{N-2}{2}} S^{N/2},\,\,B_{\mu_2}+\frac{1}{N}\mu_1^{-\frac{N-2}{2}} S^{N/2}\right\}.$$\el

\noindent {\bf Proof. } Let $N\ge 6$. For any $R>0$ small, we take $x_R\in\om$ with $\hbox{dist}(x_R, \partial\om)=3R$. Then
\be\label{eqq14}|u_{\mu_1}(x)|\le CR,\quad x\in B(x_R, 3R).\ee
Let $\psi\in C_0^{\iy}(B(0, 2))$ with $0\le\psi\le 1$
and $\psi\equiv 1$ in $B(0,1)$. Define $\vp_R(x):=1-\psi(\frac{x-x_R}{R})$, then
\be\label{cut}\vp_R(x):=\begin{cases} 0 &
\hbox{if}\quad x\in B(x_R, R), \\ 1 &
\hbox{if}\quad x\in \RN\backslash B(x_R,2R), \end {cases}\quad |\nabla \vp_R(x)|\le C/R. \ee
Define $u_R:= \vp_R u_{\mu_1}$, then by (\ref{eqq14}) and (\ref{cut}), it is easy to prove that
{\allowdisplaybreaks
\begin{gather*}
    \int_{\om}|\nabla u_R|^2\,dx\le \int_{\om}|\nabla u_{\mu_1}|^2\,dx+ CR^N;\\
    \int_{\om} |u_R|^2\,dx\ge \int_{\om}|u_{\mu_1}|^2\,dx-CR^{N+2};\\
    \int_{\om} |u_R|^{2^\ast}\,dx\ge \int_{\om}|u_{\mu_1}|^{2^\ast}\,dx-CR^{N+2^{\ast}}.
\end{gather*}
}%
Therefore, there exists $t_R>0$ independent of $\bb<0$ such that
{\allowdisplaybreaks
\begin{align*}
\max_{t>0}E_\bb(tu_R, 0)&=E_\bb(t_R u_R, 0)=\frac{1}{N}\left(\frac{\int_{\om}(|\nabla u_R|^2+\la_1 u_R^2)}{\left(\mu_1\int_{\om}
|u_R|^{2^\ast}\right)^{2/2^\ast}}\right)^{N/2}\\
&\le\frac{1}{N}\left(\frac{\int_{\om}(|\nabla u_{\mu_1}|^2+\la_1 u_{\mu_1}^2)+CR^N+CR^{N+2}}{\left(\int_{\om}\mu_1
|u_{\mu_1}|^{2^\ast}-CR^{N+2^\ast}\right)^{2/2^\ast}}\right)^{N/2}\\
&=\frac{1}{N}\left(\frac{NB_{\mu_1}+CR^N+CR^{N+2}}{\left(NB_{\mu_1}-CR^{N+2^\ast}\right)^{2/2^\ast}}\right)^{N/2}\\
&\le B_{\mu_1}+CR^N \quad \hbox{for $R>0$ small enough.}
\end{align*}
}%
Recall $U_R$ in Theorem \ref{th-5}, we have $U_R(\cdot-x_R)\cdot u_R\equiv 0$, and so $(t_R u_R, U_R(\cdot-x_R))\in\mathcal{M}_\bb$ for
all $\bb<0$. Since $N\ge 6$, one has $N>\frac{2N-4}{N-4}$. Then we see from Theorem \ref{th-5} that
{\allowdisplaybreaks
\begin{align*}
\sup_{\bb<0}B_\bb&\le E_\bb (t_R u_R, U_R(\cdot-x_R))=E_\bb (t_R u_R, 0)+E_\bb (0, U_R(\cdot-x_R))\\
&= E_\bb (t_R u_R, 0)+ J_R(U_R)\\
&\le B_{\mu_1}+CR^N+ \frac{1}{N}\mu_2^{-\frac{N-2}{2}}S^{N/2}- C_2
R^{\frac{2N-4}{N-4}}\\
&< B_{\mu_1}+ \frac{1}{N}\mu_2^{-\frac{N-2}{2}}S^{N/2} \quad \hbox{for $R>0$ small enough.}
\end{align*}
}%
By a similar argument, we also have $\sup_{\bb<0}B_\bb< B_{\mu_2}+ \frac{1}{N}\mu_1^{-\frac{N-2}{2}}S^{N/2}$.
\hfill$\square$\\

\noindent {\bf Proof of Theorem \ref{th3}. } This proof is similar to the proof of \cite[Theorem 1.4]{CZ} in case $N=4$. The novelty here is that, with the help
of Lemma \ref{lem2}, we can exclude conclusions (1)-(2) in case $N\ge 6$. Let $\beta_n<0,\,\,n\in\mathbb{N}$ satisfy $\beta_n\to-\iy$ as $n\to\iy$,
and $(u_n, v_n)$ be the positive least energy solutions of (\ref{eq2}) with $\beta=\beta_n$. By Lemma \ref{lemma1},
$E_{\beta_n}(u_n, v_n)\le A$ and so $(u_n, v_n)$ is uniformly bounded in $H$ by (\ref{eq10}).
Passing to a subsequence, we may assume that
\begin{gather*}
    u_n\rightharpoonup u_\iy,\quad v_n\rightharpoonup v_\iy\,\,\,\hbox{weakly in $H_0^1(\om)$,}
\end{gather*}
Then, by following the proof of \cite[Theorem 1.4]{CZ} in case $N=4$, we can prove that
$\int_{\om}\beta_n u_n^p v_n^p\,dx\to 0$
as $n\to\iy$, and passing to a subsequence, one of the following conclusions holds.
\begin{itemize}

\item [(1)] $u_n\to u_{\iy}$ strongly in $H_0^1(\om)$ and $v_n\rightharpoonup 0$ weakly in $H_0^1(\om)$ (so $v_n\to 0$ for almost every $x\in \om$),
where $u_\iy$ is a positive least energy
solution of
$$-\Delta u +\la_1 u=\mu_1 |u|^{2^\ast-2}u,\,\, u\in H^1_0(\om).$$
Moreover,
\be\label{eqq9}\lim_{n\to\iy}B_{\beta_n}=B_{\mu_1}+\frac{1}{N}\mu_2^{-\frac{N-2}{2}}S^{N/2}.\ee

\item [(2)] $v_n\to v_{\iy}$ strongly in $H_0^1(\om)$ and $u_n\rightharpoonup 0$ weakly in $H_0^1(\om)$ (so $u_n\to 0$ for almost every $x\in \om$),
where $v_\iy$ is a positive least energy
solution of
$$-\Delta v +\la_2 v=\mu_2 |v|^{2^\ast-2}v,\,\, v\in H^1_0(\om).$$
Moreover,
\be\label{eqq10}\lim_{n\to\iy}B_{\beta_n}=B_{\mu_2}+\frac{1}{N}\mu_1^{-\frac{N-2}{2}}S^{N/2}.\ee

\item [(3)] $(u_n, v_n)\to (u_{\iy}, v_{\iy})$ strongly in $H_0^1(\om)\times H_0^1(\om)$ and $u_\iy\cdot v_\iy= 0$ for almost $x\in\om$,  where
    $u_\iy\not\equiv 0,
v_\iy\not\equiv 0$ satisfy
\begin{align}
\label{fc1}\int_{\om}(|\nabla u_\iy|^2+\la_1 u_\iy^2)=\int_{\om}\mu_1 u_\iy^{2p},\\
\label{fc2}\int_{\om}(|\nabla v_\iy|^2+\la_2 v_\iy^2)=\int_{\om}\mu_2 v_\iy^{2p},\\
\label{least1}\lim_{n\to\iy}B_{\bb_n}= E(u_\iy, v_\iy).
\end{align}
Moreover, if $u_\iy$ and $v_\iy$ are both continuous (we will prove this later),  then $u_\iy\cdot v_\iy\equiv 0$,
 $u_\iy\in C(\overline{\om})$ is a positive least energy
solution of
$$-\Delta u +\la_1 u=\mu_1 |u|^{2^\ast-2}u,\,\, u\in H^1_0(\{u_\iy>0\}),$$
and
$v_\iy\in C(\overline{\om})$ is a positive least energy
solution of
$$-\Delta v +\la_2 v=\mu_2 |v|^{2^\ast-2}v,\,\, v\in H^1_0(\{v_\iy>0\}).$$
Furthermore, both $\{v_\iy>0\}$ and $\{u_\iy>0\}$ are connected domains, and $\{v_\iy>0\}=\om\backslash\overline{\{u_\iy>0\}}$.
\end{itemize}

Note that (\ref{eqq9})-(\ref{eqq10}) imply that one of (1) and (2) in Theorem \ref{th3} does not hold in some cases. For example, if we assume that
$-\la_1(\om)<\la_1<\la_2<0$ and $\mu_1=\mu_2$ in Theorem \ref{th3}, then
$B_{\mu_1}+\frac{1}{N}\mu_2^{-\frac{N-2}{2}}S^{N/2}<\, B_{\mu_2}+\frac{1}{N}\mu_1^{-\frac{N-2}{2}}S^{N/2},$ and so (2) in Theorem \ref{th3} does not hold,
since (\ref{eqq10}) contradicts with Lemma \ref{lemma1}.

In particular, Lemma \ref{lem2} implies that neither (1) nor (2) hold in case $N\ge 6$. That is, only (3) holds if $N\ge 6$. Therefore, the proof is complete by combining Lemma \ref{lemma7} below.\hfill$\square$\\

From the previous proof, it suffices to prove that $u_\iy$ and $v_\iy$ are continuous and $u_\iy-v_\iy$ is a least energy sign-changing solution of (\ref{fc3}). As pointed out in Remark \ref{remark3}, the following proof is completely
different from that
in \cite{CZ} for the case $N=4$.

\bl\label{lemma7} Let $(u_\iy, v_\iy)$ be in conclusion (3). Then $u_\iy-v_\iy$ is a least energy sign-changing solution of (\ref{fc3}), and $u_\iy, v_\iy$ are
both continuous.\el

\noindent {\bf Proof. } Consider the problem (\ref{fc3}).
 Its related functional is
$$P(u)=\frac{1}{2}\int_{\om}(|\nabla u|^2+\la_1 (u^+)^2+\la_2(u^-)^2)-\frac{1}{2^\ast}\int_{\om}(\mu_1 (u^+)^{2^\ast}+\mu_2 (u^-)^{2^\ast}).$$
It is standard to prove that $P\in C^1$ and its critical points are solutions of (\ref{fc3}). Define
\begin{gather*}
    J_i(u):=\int_{\om}(|\nabla u|^2+\la_i u^2-\mu_i |u|^{2^\ast}),\quad i=1, 2,\\
    \mathcal{S}:=\big\{u\in H^1_0(\om) : u^{\pm}\not\equiv 0, J_1(u^+)=0, J_2(u^-)=0\big\},\\
    m:=\inf_{u\in \mathcal{S}}P(u).
\end{gather*}
Then any sign-changing solutions of (\ref{fc3}) belong to $\mathcal{S}$. By (\ref{fc1})-(\ref{fc2}), we have $u_\iy-v_\iy\in \mathcal{S}$ and so $m\le
P(u_\iy-v_\iy)=E(u_\iy, v_\iy)$.
For any $u\in \mathcal{S}$, we have $(u^+, u^-)\in \mathcal{M}_{\bb}$ for all $\bb$. Then by (\ref{least1}) we see that
\begin{align*}
P(u_\iy-v_\iy)=E(u_\iy, v_\iy)=\lim_{n\to\iy}B_{\bb_n}\le E(u^+, u^-)=P(u),\quad\forall\, u\in \mathcal{S},
\end{align*}
and so $P(u_\iy-v_\iy)\le m$. Combining these with Lemma \ref{lem2}, we obtain
\be\label{fc4}P(u_\iy-v_\iy)=m<\min\left\{B_{\mu_1}+\frac{1}{N}\mu_2^{-\frac{N-2}{2}} S^{\frac{N}{2}},\,\,B_{\mu_2}+\frac{1}{N}\mu_1^{-\frac{N-2}{2}} S^{\frac{N}{2}}\right\}.\ee

{\bf Step 1.} We show that $P'(u_\iy-v_\iy)=0$, and so $u_\iy-v_\iy$ is a least energy sign-changing solution of (\ref{fc3}).

Thanks to (\ref{fc4}), the following argument is standard (see \cite{LWZ, SWT} for example), and we give the details here for completeness.

Assume that $u_\iy-v_\iy$ is not a critical point of $P$, then there exists $\phi\in C^\iy_0(\om)$ such that $P'(u_\iy-v_\iy)\phi\le -1$. Then there exists
$0<\e_0<1/10$, such that for $|t-1|\le \e_0$, $|s-1|\le \e_0$, $|\sg|\le \e_0$, there holds
$$P'(tu_\iy-s v_\iy+\sg \phi)\phi\le -\frac{1}{2}.$$
Consider a function $0\le\eta\le 1$ defined for $(t, s)\in T=[\frac{1}{2}, \frac{3}{2}]\times[\frac{1}{2}, \frac{3}{2}]$, such that
\begin{align*}
&\eta(t, s)=1, \quad \hbox{for $|t-1|\le\frac{\e_0}{2}, \,\,|s-1|\le \frac{\e_0}{2}$},\\
&\eta(t, s)=0, \quad \hbox{for $|t-1|\ge\e_0$ or $|s-1|\ge\e_0$}.
\end{align*}
Then for $|t-1|\le\e_0, |s-1|\le \e_0$, we have
\begin{align*}
&P(tu_\iy-sv_\iy +\e_0\eta(t, s)\phi)\\
=&P(tu_\iy-sv_\iy)+\int_0^1P'(tu_\iy-sv_\iy +\theta\e_0\eta(t, s)\phi)[\e_0 \eta(t, s)\phi]d\theta\\
\le& P(tu_\iy-sv_\iy)-\frac{1}{2}\e_0 \eta(t, s).
\end{align*}
Note that
$$\sup_{t, s>0} P(t u_\iy-s v_\iy)=P(u_\iy-v_\iy)=m,$$
and for $|t-1|\ge\frac{\e_0}{2}$ or $|s-1|\ge\frac{\e_0}{2}$, there exists $0<\dd<\frac{\e_0}{2}$ such that
$$P(t u_\iy-s v_\iy)\le m-\dd.$$
We have, for $|t-1|\le\frac{\e_0}{2}, |s-1|\le \frac{\e_0}{2}$, that
$$P(tu_\iy-sv_\iy +\e_0\eta(t, s)\phi)\le m-\frac{\e_0}{2};$$
for $\frac{\e_0}{2}\le|t-1|\le\e_0, |s-1|\le \e_0$ or $\frac{\e_0}{2}\le|s-1|\le\e_0, |t-1|\le \e_0$,
$$P(tu_\iy-sv_\iy +\e_0\eta(t, s)\phi)\le P(tu_\iy-sv_\iy)\le m-\dd;$$
for $|t-1|\ge \e_0$ or $|s-1|\ge \e_0$,
$$P(tu_\iy-sv_\iy +\e_0\eta(t, s)\phi)= P(tu_\iy-sv_\iy)\le m-\dd.$$
So
\be\label{fc5}\sup_{(t, s)\in T}P(t u_\iy-s v_\iy+\e_0\eta(t, s)\phi)\le m-\dd.\ee
On the other hand, for $\e\in [0, \e_0]$, let $h_\e : T\to H^1_0(\om)$ by $h_\e(t, s)=t u_\iy-s v_\iy+\e \eta(t, s)\phi$, and
$H_\e: T\to \R^2$ by
$$H_\e(t, s)=(J_1(h_\e(t, s)^+), J_2(h_\e(t, s)^-)).$$
Note that for any $(t, s)\in \partial T$, we have $\eta(t, s)=0$ and so $h_\e(t, s)\equiv h_0(t, s)=tu_\iy-sv_\iy$. Moreover,
$H_0(t, s)=(J_1(t u_\iy),\,\, J_2(s v_\iy))$.
Then it is easy to see that
$$\hbox{deg}(H_{\e_0}(t, s), T, (0, 0))=\hbox{deg}(H_{0}(t, s), T, (0, 0))=1,$$
that is, there exists $(t_0, s_0)\in T$ such that $h_{\e_0}(t_0, s_0)\in \mathcal{S}$, which is a contradiction with (\ref{fc5}).

{\bf Step 2.} We show that $u_\iy$ and $v_\iy$ are continuous.

By Step 1, $u_\iy-v_\iy$ is a nontrivial solution of (\ref{fc3}). Then by a Brezis-Kato argument (see \cite{BK}), we see that
$u_\iy-v_\iy\in L^q(\om), \,\forall\, q\ge 2$. In particular,
$$\mu_1 u_\iy^{2^\ast-1}-\mu_2 v_\iy^{2^\ast-1}-\la_1 u_\iy+\la_2 v_\iy\in L^q(\om), \,\,\forall \, q > N.$$
Then by elliptic regularity theory, $u_\iy-v_\iy\in W^{2, q}(\om)$ with $q>N$. By Sobolev embedding, we have $u_\iy-v_\iy\in C(\overline{\om})$. Since
$u_\iy=(u_\iy-v_\iy)^+$ and $v_\iy=(u_\iy-v_\iy)^-$, we see that $u_\iy$ and $v_\iy$ are both continuous. This completes the proof and so completes the proof of
Theorem \ref{th3}.
\hfill$\square$\\

\noindent {\bf Proof of Theorem \ref{sign-changing}. } Let $N\ge 6$. Actually, by Theorem \ref{th3} and Lemma \ref{lemma7}, we have proved that the problem
(\ref{fc3}) has a least energy sign-changing solution $u_\iy-v_\iy$. Obviously,
Theorem \ref{sign-changing} is a direct corollary by letting $\la_1=\la_2$ and $\mu_1=\mu_2$, and (\ref{signchange}) follows directly from (\ref{fc4}). \hfill$\square$

\end{document}